\documentclass[11pt, twoside]{article}
\usepackage{mathrsfs}
\usepackage{amsfonts}
\usepackage{setspace}
\usepackage{amsmath,amssymb}
\usepackage{multicol}
\usepackage{color}
\usepackage{cite}
\usepackage{graphicx}
\newtheorem{theorem}{Theorem}[section]
\newtheorem{lemma}[theorem]{Lemma}

\newtheorem{definition}[theorem]{Definition}
\newtheorem{remark}[theorem]{Remark}

\numberwithin{equation}{section}
\newenvironment{proof}[1][Proof]{\noindent\textbf{#1. }}{\hfill $\Box$}
\allowdisplaybreaks

\newcommand{\beq}{\begin{equation}}
\newcommand{\eeq}{\end{equation}}
\newcommand{\eqnref}[1]{(\ref {#1})}
\makeatletter
\begin{document}

\author{Youjun Deng\footnote{School of Mathematics and Statistics, Central South University, Changsha, Hunan, China. Email: youjundeng@csu.edu.cn; dengyijun\_001@163.com}\,,\ Hongyu Liu\footnote{Department of Mathematics, City University of Hong Kong, Kowloon, Hong Kong SAR, China. Email: hongyu.liuip@gmail.com; hongyliu@cityu.edu.hk}\ \, and\ Guang-Hui Zheng\footnote{School of Mathematics, Hunan University, Changsha, Hunan, China. Email: zhenggh2012@hnu.edu.cn}
}

\date{}

\title{\textbf{\Large Mathematical analysis of plasmon resonances for curved nanorods}}
\maketitle
\begin{abstract}
We investigate plasmon resonances for curved nanorods which present anisotropic geometries. We analyze quantitative properties of the plasmon resonance and its relationship to the metamaterial configurations and the anisotropic geometries of the nanorods. Based on delicate and subtle asymptotic and spectral analysis of the layer potential operators, particularly the Neumann-Poincar\'e operators, associated with anisotropic geometries, we derive sharp asymptotic formulae of the corresponding scattering field in the quasi-static regime. By carefully analyzing the asymptotic formulae, we establish sharp conditions that can ensure the occurrence of the plasmonic resonance. The resonance conditions couple the metamaterial parameters, the wave frequency and the nanorod geometry in an intricate but elegant manner. We provide thorough resonance analysis by studying the wave fields both inside and outside the nanorod. Furthermore, our quantitative analysis indicates that different parts of the nanorod induce varying degrees of resonance. Specifically, the resonant strength at the two end-parts of the curved nanorod is more outstanding than that of the facade-part of the nanorod. This paper presents the first theoretical study on plasmon resonances for nanostructures within anisotropic geometries.

\noindent\textbf{Keywords}: Plasmon resonance; curved nanorod; anisotropic geometry; asymptotic analysis; spectral analysis; Neumann-Poincar\'{e} operator

\noindent\textbf{Mathematics Subject Classification:}~~35Q60, 35J05, 31B10, 35R30, 78A40

\end{abstract}

\tableofcontents
\section{Introduction}

\subsection{Physical background}

Localized surface plasmon is the resonant oscillation of conducting electrons at the interface between negative and positive permittivity material stimulated by incident light.
Owing to their unique optical or electromagnetic properties, plasmonic materials have been applied in various scientific fields, such as enhancing the brightness of light, confining strong electromagnetic fields, medical therapy, invisibility cloaking, biomedical imaging \cite{ACKLM5,CKKL7,SC1,KLSW8,JLEE3,BGQ2,LLL9,LL15,LLL16,WN10,BS11,ADM6,Z22} and so on. The plasmon technology is revolutionizing many industrial applications.

Metallic nanostructures including nanoparticles, nanorods and photonic crystal form the basis for constructing various plasmonic devices. A metallic nanoparticle is the simplest nanostructure, which possesses a small size uniformly in all dimensions. Recently, there are extensive and intensive studies on mathematically characterizing the plasmon resonances associated with nanoparticles; see \cite{ACKLM5,ACL20,CKKL7,SC1,KLSW8,AKL13,AMRZ14,FDL15,LLL9,LL15,LLL16,WN10,BS11,ADM6, Z22} and the references cited therein. The metallic nanorod is another important nanostructure, which has a long aspect ratio and possesses different size scales in different dimensions. That is, the nanorods present anisotropic geometries. Metallic nanorods have been widely used in real applications including semiconductor materials, microelectromechanical systems, food packaging, catalysis, energy storage, biomedicine and cloaking \cite{LWJC41,JGM40,HDA42,VTBH44,D43}. In particular, some nanorods (such as $C_eO_2$) display higher catalytic activity compared to the nanoparticles, which could potentially increase their usage \cite{LWJC41}. However, to our best knowledge, there is little mathematical study in the literature on theoretically characterizing the plasmon resonances associated with nonorods. Most of the existing studies in the aforementioned physics literature are concerned with experimental/practical applications associated with specific nanorod geometries such as a long and straight cylinder or a long and thin ellipsoid. In this paper, we aim to present a comprehensive mathematical analysis on the plasmon resonances associated with general nanorod geometries, dubbed as the curved nanorods. We establish the sharp resonance conditions and accurately characterize the quantitative behaviours of the resonant fields. Our study connects to the spectral properties of the layer potential operators, particularly the Neumann-Poincar\'e operators, associated with the anisotropic nanorod geometries. Compared to the related studies for nanoparticles (isotropic geometries), the corresponding asymptotic and spectral analysis associated with nanorods (anisotropic geometries) are radically more challenging. The study presented in this work not only supplies a fundamental basis for many real applications of plasmon nanorod structures, but also opens up a field of research with many potential developments. Our analysis is mainly conducted for the Helmholtz system and the mathematical strategies developed pave the way for extensions to the other physical systems, in particular the full Maxwell system.

\subsection{Mathematical formulation and setup}

We present the mathematical setup of our study by first introducing the anisotropic nanorod geometry. We shall adopt the notations in \cite{DLU7, DLU72} where anisotropic geometries were introduced in a different context. Let $\Gamma_0$ be a smooth simple and non-closed curve in $\mathbb{R}^3$, and the two endpoints of $\Gamma_0$ are $P_0$ and $Q_0$. Let $r\in\mathbb{R}_+$.
Denote by $N(x)$ the normal plane of the curve $\Gamma_0$ at $x\in\Gamma_0$. We note that $N(P_0)$ and $N(Q_0)$ are, respectively, defined by the left and right limits along $\Gamma_0$. For any $x\in\Gamma_0$, we let $\mathscr{S}_r(x)$ denote the disk lying on $N(x)$, centered at $x$ and of radius $r$. It is assumed that there exists $r_0\in\mathbb{R}_+$ such that when $r\leq r_0$, $\mathscr{S}_r(x)$ intersects $\Gamma_0$ only at $x$. We start with a thin structure $D_r^f$ given by
\begin{equation}\label{2.1}
D_r^f:=\mathscr{S}_r(x)\times\Gamma_0(x),\ x\in\overline{\Gamma}_0,
\end{equation}
where we identify $\Gamma_0$ with its parametric representation $\Gamma_0(x)$. Clearly, the facade of $D_r^f$, denoted by $S_r^f$ and parallel to $\Gamma_0$, is given by
\begin{equation}\label{2.2}
S_r^f:=\{x+r\cdot \mathbf{n}(x); x\in \Gamma_0, \mathbf{n}(x)\in N(x) \cap \mathbb{S}^2\},
\end{equation}
and the two end-surfaces of $D_r^f$ are the two disks $\mathscr{S}_r(P_0)$ and $\mathscr{S}_r(Q_0)$. Let $D_{r_0}^a$ and $D_{r_0}^b$ be two simply connected sets with $\partial D_{r_0}^a=S_{r_0}^a\cup \mathscr{S}_{r_0}(P_0)$ and $\partial D_{r_0}^b=S_{r_0}^b\cup\mathscr{S}_{r_0}(Q_0)$. It is assumed that $S_{r_0}:=S^f_{r_0}\cup S^b_{r_0} \cup S^a_{r_0}$ is a smooth boundary of the domain $D_{r_0}:=D_{r_0}^a\cup D_{r_0}^f\cup D_{r_0}^b$. For $0<r<r_0$, we set
\[
D_r^a:=\frac{r}{r_0}(D_{r_0}^a-P_0)+P_0=\left\{\frac{r}{r_0}\cdot(x-P_0)+P_0; x\in D_{r_0}^a\right\},
\]
and similarly, $D_r^b:={r}/{r_0}\cdot(D_{r_0}^b-Q_0)+Q_0$. Let $S_r^a$ and $S_r^b$, respectively, denote the boundaries of $D_r^a$ and $D_r^b$ excluding $\mathscr{S}_r^a$ and $\mathscr{S}_r^b$. Now, we set $D_r:=D_r^a\cup D_r^f\cup D_r^b$, and $S_r:=S^f_r\cup S^b_r \cup S^a_r=\partial D_r$. $D_r$ represents the geometry of a curved nanorod in our study with $D_r^{a, b}$ signifying the two end-parts and $D_r^f$ signifying the facade-part.

According to our earlier description, it is obvious that for $0<r\leq r_0$, $D_r$ is a simply connected set with a smooth boundary $S_r$,  and $D_{r_1}\Subset D_{r_2}$ if $0\leq r_1<r_2\leq r_0$. Moreover,  $D_r$ degenerates to $\Gamma_0$ if one takes $r=0$. Without loss of generality, we assume that $r_0\equiv 1$. In what follows, we let $\delta\in\mathbb{R}_+$ be the size parameter and let
\begin{equation}\label{2.3}
S_\delta:=S^f_\delta \cup S^b_\delta \cup S^a_\delta,
\end{equation}
denote the boundary surface of $D_\delta$. In order to ease the exposition, we drop the dependence on $r$ if one takes $r=1$. For example, $D$ and $S$ denote, respectively, $D_r$ and $S_r$ with $r=1$. It is emphasized that in all of our subsequent argument, $D$ can always be replaced by $D_{\tau_0}$ with $0<\tau_0\leq r_0$ being a fixed number. Hence, we indeed shall not lose any generality of our study by assuming that $r_0\equiv 1$.
Finally, we would like to note that
a particular case is to take $\Gamma_0$ to be a straight line-segment and,
$S^{b}_r$ and $S^{a}_r$ to be two semi-spheres of radius $r$ and centered at $P_0$ and $Q_0$ respectively. Hence, though we term $D_\delta$ as a curved nanorod, it actually also includes as a special case the ``straight" nanorod as considered in the physics literature.

Next we introduce a blowup transformation which maps $y\in\overline{D_\delta}$ to $\tilde y\in \overline{D}$ as follows
\begin{equation}\label{2.4}
A(y)=\tilde y:=\frac{1}{\delta}(y-z_y)+z_y,\quad y\in D_\delta^f,
\end{equation}
whereas
\begin{equation}\label{2.5}
A(y)=\tilde y:= \left\{
\begin{array}{ll}
\frac{y-P_0}{\delta}+P_0, & y \in D_\delta^a, \\
\frac{y-Q_0}{\delta}+Q_0, & y \in D_\delta^b,
\end{array}
\right.
\end{equation}
which will be used to analyze the asymptotic expansions of some layer potential operators.

In this paper, we consider a long and thin nanorod occupying a bounded and simply connected domain $D_\delta\subset\mathbb{R}^3$ as described above, whose boundary $\partial {D_\delta}\ (\equiv S_\delta)$ is $\mathcal{C}^{1,\gamma}$ for some $\gamma\in(0,1)$. We present the mathematical description of the electromagnetic (EM) scattering from the nanorod $D_\delta$. In principle, the propagation of light in the nanostructures is described by the Maxwell equations. Due to technical reasons, we shall mainly consider in this paper the Helmholtz equation, which in 2D describes the transverse EM propagation and in 3D the acoustic wave propagation (cf. \cite{AKL13,LL15}). In order to avoid repeating the discussions, we present the results mainly for the 3D case and the extension to the 2D case should be clear and can be straightforwardly done. As remarked earlier, the extension to the full Maxwell system will be presented in our forthcoming work. In order to ease the exposition, we stick to the physical terminologies associated with the electromagnetic scattering and the physical interpretation of our results for the acoustic scattering can be easily given.

 The material properties of the nanorod $D_\delta$ are characterized by the electric permittivity $\varepsilon_c$ and the magnetic permeability $\mu_c$, while the homogeneous medium in $\mathbb{R}^3\setminus\overline{{D_\delta}}$ is characterized by electric permittivity $\varepsilon_m$ and the magnetic permeability $\mu_m$. We shall be mainly concerned with the time-harmonic scattering and let $\omega\in\mathbb{R}_+$ signify the angular frequency of the wave.  Let $\Re\varepsilon_c<0$, $\Im\varepsilon_c>0$, $\Re\mu_c<0$, $\Im\mu_c>0$ be constants. Define the wavenumbers to be
\begin{equation}\label{2.6}
k_c=\omega\sqrt{\varepsilon_c\mu_c},\ \ k_m=\omega\sqrt{\varepsilon_m\mu_m}.
\end{equation}
Set
\begin{equation}\label{2.7}
\varepsilon_{D_\delta}=\varepsilon_c\chi({D_\delta})+\varepsilon_m\chi(\mathbb{R}^3\setminus\overline{{D_\delta}}),\ \ \mu_{D_\delta}=\mu_c\chi({D_\delta})+\mu_m\chi(\mathbb{R}^3\setminus\overline{{D_\delta}}),
\end{equation}
where and also in what follows, $\chi$ signifies the characteristic function. It is pointed out that nano-metal materials with the electric permittivity and magnetic permeability satisfying $\Re\varepsilon_c<0$, $\Re\mu_c<0$ are called double negative materials, which show several unusual properties, such as the counter directance between the group velocity and the phase vector, negative index of refraction and the reverse Doppler and Cherenkov effects \cite{ZSK20,V21,SPW22}. It is emphasized that the double negative materials inside the nanorod are the key ingredient accounting for the plasmon resonance in our study. The positive imaginary parts of electric permittivity and magnetic permeability signify the dissipation of the plasmonic nanostructures. For the background medium in $\mathbb{R}^3\backslash\overline{D_\delta}$, we assume that $\varepsilon_m$, $\mu_m$ are real and strictly positive.

Let $u^i=e^{\mathrm{i}k_m d\cdot x}$, $\mathrm{i}:=\sqrt{-1}$, be the time-harmonic incident plane wave. Here $d$ is a unit vector which represents the incident direction. The wave scattering due to the impingement of the incident field $u^i$ on the nanorod $D_\delta$ is described by the following Helmholtz system,
\begin{align}
\begin{cases}\label{2.8}
\nabla\cdot\frac{1}{\varepsilon_{D_\delta}}\nabla u+\omega^2\mu_{D_\delta} u=0 &\text{in}\ \mathbb{R}^3\setminus\partial {D_\delta},\\
u|_+=u|_- &\text{on}\ \partial {D_\delta},\\
\frac{1}{\varepsilon_m}\frac{\partial
u}{\partial\nu}\big|_+=\frac{1}{\varepsilon_c}\frac{\partial
u}{\partial\nu}\big|_- &\text{on}\ \partial {D_\delta}\\
u^s:=u-u^i &\text{satisfies the Sommerfeld radiation condition}.
\end{cases}
\end{align}
By the Sommerfeld radiation condition, we mean that the scattered wave $u^s$ satisfies
\begin{equation}\label{2.9}
\frac{\partial u^s}{\partial |x|}-\mathrm{i}k_m u^s=\mathcal{O}(|x|^{-2})\quad\mbox{as $|x|\rightarrow +\infty$},
\end{equation}
which holds uniformly in the angular variable $\hat x=x/|x|$. The Sommerfeld radiation condition characterises the outgoing nature of the scattered field.

Our study of the plasmon resonance associated with the Helmholtz system \eqref{2.8} heavily relies on the layer-potential techniques. To that end, we introduce the boundary layer potential operators for the analysis of the solution to \eqnref{2.8}, and we also refer to \cite{AGJKLSW28,CK26} for more relevant discussions on layer potential techniques. Let $G^{k}(x,y)$ be the 3D Green's function for the PDO $\Delta+k^2$, which is given by
\begin{equation}\label{2.10}
G^{k}(x,y)=-\frac{e^{\mathrm{i}k|x-y|}}{4\pi|x-y|}.
\end{equation}
Let $\Sigma$ be a bounded domain with a $C^{1,\gamma}$ boundary $\partial\Sigma$, which could be $D_\delta$ or $D$ in our subsequent study.
The single layer potential for the Helmholtz equation is defined by
\begin{equation}\label{2.11}
\mathcal {S}_\Sigma^{k}[\phi](x)=\int_{\partial \Sigma}G^{k}(x,y)\phi(y)d\sigma(y),\ \ \ x\in \mathbb{R}^3,
\end{equation}
where $\phi\in H^{-\frac{1}{2}}(\partial \Sigma)$ signifies a boundary density function.
The following jump formula holds
\begin{equation}\label{2.12}
\frac{\partial(\mathcal{S}_\Sigma^{k}[\phi])}{\partial\nu}\bigg|_{\pm}(x)=\left(\pm\frac{1}{2}\mathcal
{I}+(\mathcal {K}_\Sigma^{k})^*\right)[\phi](x),\ \ \ a.e.\ x\in\partial \Sigma,
\end{equation}
where
\begin{equation}\label{2.13}
(\mathcal{K}_\Sigma^{k})^*[\phi](x)= \mbox{p.v.}\quad \int_{\partial \Sigma}\frac{\partial G^{k}(x,y)}{\partial\nu(x)}\phi(y)d\sigma(y),
\end{equation}
is known as the Neumann-Poincar\'e operator. Here p.v. stands for the Cauchy principle value. In what follows, we let $\mathcal{S}_\Sigma$ and $\mathcal{K}^*_\Sigma$ respectively denote the operators $\mathcal{S}^k_\Sigma$ and $(\mathcal{K}^k_\Sigma)^*$ by formally taking $k=0$, which are known as the static single-layer and Neumann-Poincar\'e operators, respectively.

By using the single layer potential (\ref{2.11}) and the jump formula (\ref{2.12}), one has the following integral representation for the solution to (\ref{2.8}):
\begin{align}\label{2.14}
u(x)=
\begin{cases}
\mathcal {S}_{D_\delta}^{k_c}[\phi](x), &x\in {D_\delta},\\
u^i(x)+\mathcal {S}_{D_\delta}^{k_m}[\psi](x), &x\in
\mathbb{R}^3\setminus\overline{{D_\delta}},
\end{cases}
\end{align}
where $(\phi,\psi)\in H^{-\frac{1}{2}}(\partial {D_\delta})\times H^{-\frac{1}{2}}(\partial {D_\delta})$ satisfy the following integral system
\begin{equation}\label{2.15}
\begin{cases}
\mathcal {S}_{D_\delta}^{k_c}[\phi]-\mathcal {S}_{D_\delta}^{k_m}[\psi]=u^i &\text{on}\ \partial {D_\delta},\medskip\\
\frac{1}{\varepsilon_c}\left(-\frac{1}{2}\mathcal {I}+(\mathcal
{K}_{D_\delta}^{k_c})^*\right)[\phi]
-\frac{1}{\varepsilon_m}\left(\frac{1}{2}\mathcal {I}+(\mathcal
{K}_{D_\delta}^{k_m})^*\right)[\psi] =\frac{1}{\varepsilon_m}\frac{\partial u^i}{\partial\nu} &\text{on}\ \partial {D_\delta}.
\end{cases}
\end{equation}

We focus on our analysis in the quasi-static regime, namely $\omega\cdot \mathrm{diam}(D_\delta)\ll1$. By a standard scaling argument and without loss of generality, we could assume throughout the rest of the paper that $\mathrm{diam}(D_\delta)\sim 1$ and $\omega\ll 1$. It is emphasised that $\delta\ll 1$ is also an asymptotic parameter and hence $\mathrm{diam}(D_\delta)\sim \mathrm{diam}(D_\delta^f)$. In fact, $\delta$ is the key parameter that characterises the anisotropy of the geometry of the nanorod $D_\delta$ and it shall be related to $\omega$ in what follows. This should become more evident in our subsequent discussion.

It is known that for $\omega$ small enough, $\mathcal {S}_{D_\delta}^{k}$ is invertible (cf. \cite{DLU7}). Therefore, by using the first equation in (\ref{2.15}), one can directly obtain that
\begin{equation}\label{2.16}
\phi=\left(\mathcal {S}_{D_\delta}^{k_c}\right)^{-1}\left(\mathcal {S}_{D_\delta}^{k_m}[\psi]+u^i\right).
\end{equation}
Then, from the second equation in (\ref{2.15}), we have that
\begin{equation}\label{2.17}
\mathcal {A}_{D_\delta}(\omega)[\psi]=f,
\end{equation}
where
\begin{align}\label{2.18}
\mathcal {A}_{D_\delta}(\omega)=&\frac{1}{\varepsilon_m}\left(\frac{1}{2}\mathcal
{I}+(\mathcal {K}_{D_\delta}^{k_m})^*\right)
+\frac{1}{\varepsilon_c}\left(\frac{1}{2}\mathcal {I}-(\mathcal
{K}_{D_\delta}^{k_c})^*\right)(\mathcal {S}_{D_\delta}^{k_c})^{-1}\mathcal
{S}_{D_\delta}^{k_m},\\
f=&-\frac{1}{\varepsilon_m}\frac{\partial u^i}{\partial\nu}-\frac{1}{\varepsilon_c}\left(\frac{1}{2}\mathcal {I}-(\mathcal
{K}_{D_\delta}^{k_c})^*\right)(\mathcal {S}_{D_\delta}^{k_c})^{-1}[u^i].\label{2.19}
\end{align}
Clearly,
\begin{align}
\mathcal {A}_{D_\delta}(0)=\mathcal {A}_{D_\delta,0}=&\frac{1}{\varepsilon_m}\left(\frac{1}{2}\mathcal {I}+\mathcal
{K}_{D_\delta}^*\right) +\frac{1}{\varepsilon_c}\left(\frac{1}{2}\mathcal
{I}-\mathcal
{K}_{D_\delta}^*\right)\nonumber\\
=&\frac{1}{2}\left(\frac{1}{\varepsilon_m}+\frac{1}{\varepsilon_c}\right)\mathcal
{I}+\left(\frac{1}{\varepsilon_m}-\frac{1}{\varepsilon_c}\right)\mathcal
{K}_{D_\delta}^*.\label{2.20}
\end{align}

Similarly, for $\omega\ll1$, we can also deduce the following operator equation of the density $\phi$:
\begin{equation}\label{2.21}
\widetilde{\mathcal {A}}_{D_\delta}(\omega)[\phi]=\tilde{f},
\end{equation}
where
\begin{align}\label{2.22}
\widetilde{\mathcal {A}}_{D_\delta}(\omega)=&\frac{1}{\varepsilon_c}\left(\frac{1}{2}\mathcal
{I}-(\mathcal {K}_{D_\delta}^{k_c})^*\right)
+\frac{1}{\varepsilon_m}\left(\frac{1}{2}\mathcal {I}+(\mathcal
{K}_{D_\delta}^{k_m})^*\right)(\mathcal {S}_{D_\delta}^{k_m})^{-1}\mathcal
{S}_{D_\delta}^{k_c},\\
\tilde{f}=&-\frac{1}{\varepsilon_m}\frac{\partial u^i}{\partial\nu}+\frac{1}{\varepsilon_m}\left(\frac{1}{2}\mathcal {I}+(\mathcal
{K}_{D_\delta}^{k_m})^*\right)(\mathcal {S}_{D_\delta}^{k_m})^{-1}[u^i].\label{2.23}
\end{align}
Notice that $\widetilde{\mathcal {A}}_{D_\delta}(0)=\widetilde{\mathcal {A}}_{D_\delta,0}=\mathcal {A}_{D_\delta}(0)$.


Finally, we introduce the formal definition of the plasmon resonance for our subsequent study.

\begin{definition}\label{depr}
Consider the wave scattering system \eqref{2.8} associated with the nanorod $D_\delta$, where the material configuration is described in \eqref{2.6} and \eqref{2.7}.
\begin{enumerate}
\item[(i)] Plasmon resonance occurs if the following condition is fulfilled:
\begin{align}\label{prdf01}
\left\|\nabla u^s\right\|_{L^2(\mathbb{R}^3\setminus\overline{D_\delta})}\gg 1.
\end{align}

\item[(ii)] If the internal energy of wave field inside the nanorod significantly increases, i.e.,
\begin{align}\label{prdf02}
\left\|\nabla u\right\|_{L^2(D_\delta)}\gg 1,
\end{align}
then we also say that plasmon resonance occurs.

\item[(iii)] Define the electric energy of the nanorod as follows (cf. \cite{BW39}):
\begin{align}\label{ee}
\mathcal{E}=\frac{1}{8\pi}\Im(\varepsilon_c)\left\|\nabla u\right\|_{L^2(D_\delta)}^2.
\end{align}
We say that plasmon resonance occurs with respect to the electric energy if $\mathcal{E}\gg 1$.
\end{enumerate}
\end{definition}

\begin{remark}
According to \eqref{prdf01}, \eqref{prdf02} or \eqref{ee}, we see that when plasmon resonance occurs, the resonant wave field exhibits highly oscillatory behaviours that cause the blowup of the resonant energy in different senses. The high oscillation is a hallmark of the plasmon resonance and moreover the high oscillation is mainly confined within the vicinity of the nanostructure, which is the fundamental basis for many plasmonic technologies.
\end{remark}

\begin{remark}
We would like to emphasize that compared to the other definitions in the mathematical literature on the plasmon resonance, say e.g. Definition 1 in \cite{AMRZ14}, our definition of the plasmon resonance is a more refined one. Indeed, in Definition~\ref{depr}, we characterize the resonance with respect to the wave fields both outside and inside the nanorod, as well as the electric energy. It is pointed out that in Definition \ref{depr}, (i) and (ii) are essentially equivalent; see Theorems \ref{thm4.3} and \ref{thminter}. (iii) means that stronger resonant behaviours occur since one usually has $\Im(\varepsilon_c)\ll 1$. Nevertheless, through our analysis, we can establish conditions such that (i)-(iii) can occur concurrently; see Remark \ref{remee} in what follows.
\end{remark}


\section{Asymptotic and quantitative analysis of the scattering field}

In this section, we conduct asymptotic analysis of the scattering system \eqref{2.8} associated with the nanorod $D_\delta$. We shall derive several asymptotic formulas of the wave field $u^s$ with respect to $\omega\ll 1$ and $\delta\ll 1$, which are crucial to our subsequent analysis of the resonant behaviours of the field, in particular the anisotropy of the resonant behaviours that is related to the anisotropic geometry of the nanorod.

\subsection{Asymptotics of layer potential operators}\label{sec2}

In this section, we derive several asymptotic expansion formulas of the layer potential operators with respect to the asymptotic wavenumber $k$ as well as the anisotropic size parameter $\delta$ of the nanorod. Those asymptotic results pave the way for the asymptotic analysis of the scattering system \eqref{2.8}. It is pointed out that some of the asymptotic results can be found in \cite{AMRZ14,DLU7}, which we still include for the self-containedness, and some are new from the current setup of study whose proofs shall be given.

First, it is recalled the layer potential operators $\mathcal{S}^k_\Sigma, (\mathcal{K}^k_\Sigma)^*$ and $\mathcal{S}_\Sigma, \mathcal{K}^*_\Sigma$ introduced earlier. For the subsequent use, we introduce the function space $\mathcal {H}^*(\partial \Sigma)$ which is a Hilbert space equipped
with the following inner product
\begin{align}\label{3.1}
\langle u,v\rangle_{\mathcal {H}^*(\partial \Sigma)}=-(u,\mathcal
{S}_{\Sigma}[v])_{-\frac{1}{2},\frac{1}{2}},
\end{align}
where $(\cdot,\cdot)_{-\frac{1}{2},\frac{1}{2}}$ is the duality pairing between the Sobolev spaces $H^{-\frac{1}{2}}(\partial \Sigma)$ and $H^{\frac{1}{2}}(\partial \Sigma)$. It is noted that $\|\cdot\|_{\mathcal{H}^*(\partial\Sigma)}$ is equivalent to $\|\cdot\|_{H^{-1/2}(\partial\Sigma)}$ (cf. \cite{ACKLM5}).

The following lemma collects the necessary asymptotic results of the layer potential operators with respect to $k\ll 1$, whose proof can be found in \cite{AMRZ14}.
\begin{lemma}\label{lem:a1}
For $k\ll 1$, the following asymptotic results hold.
\begin{enumerate}
\item[(i)] The single layer potential operator can be expanded as follows:
\begin{equation}\label{8.1}
\mathcal {S}_D^{k}=\mathcal{S}_D+\sum_{j=1}^\infty k^{j}\mathcal{S}_{D,j},
\end{equation}
where
\begin{equation}\label{sing01}
\mathcal{S}_{D,j}[\psi](x)=-\frac{\mathrm{i}^j}{4\pi j!}\int_{\partial D}|x-y|^{j-1}\psi(y)d\sigma(y).
\end{equation}
Moreover, the norms $\|\mathcal{S}_{D,j}\|_{\mathcal{L}(\mathcal {H}^*(\partial D),\mathcal
{H}^*(\partial D))}$ are uniformly bounded with respect to $j\geq1$, and the series in (\ref{8.1}) is convergent in
{$\mathcal {L}(\mathcal {H}^*(\partial D),\mathcal{H}^*(\partial D))$. }

\item[(ii)] It holds that
\begin{equation}\label{8.3}
(\mathcal {S}_D^{k})^{-1}=\mathcal{S}_D^{-1}+k\mathcal{B}_{D,1}+k^2\mathcal{B}_{D,2}+\cdots,
\end{equation}
where
\begin{align}
\mathcal{B}_{D,1}=-\mathcal{S}_D^{-1}\mathcal{S}_{D,1}\mathcal{S}_D^{-1},\quad
\mathcal{B}_{D,2}=-\mathcal{S}_D^{-1}\mathcal{S}_{D,2}\mathcal{S}_D^{-1}+\mathcal{S}_D^{-1}\mathcal{S}_{D,1}\mathcal{S}_D^{-1}\mathcal{S}_{D,1}\mathcal{S}_D^{-1}.
\end{align}
Furthermore, the series in (\ref{8.3}) is convergent in
$\mathcal{L}(\mathcal {H}^*(\partial D),\mathcal{H}^*(\partial D))$.

\item[(iii)] The Neumann-Poincar\'e operator $(\mathcal{K}_D^k)^*$ can be expanded as follows:
\begin{equation}\label{8.6}
(\mathcal {K}_D^{k})^*=\mathcal{K}_D^*+\sum_{j=1}^\infty k^{j}\mathcal{K}_{D,j},
\end{equation}
where
\begin{equation}\label{doub01}
\mathcal{K}_{D,j}[\psi](x)=-\frac{\mathrm{i}^j(j-1)}{4\pi j!}\int_{\partial D}|x-y|^{j-3}( x-y,\nu(x))\psi(y)d\sigma(y)
\end{equation}
Furthermore, the norms $\|\mathcal{S}_{D,j}\|_{\mathcal{L}(\mathcal {H}^*(\partial D),\mathcal
{H}^*(\partial D))}$ are uniformly bounded with respect to $j$, and the series in (\ref{8.6}) is convergent in
 $\mathcal {L}(\mathcal {H}^*(\partial D),\mathcal{H}^*(\partial D))$.
\end{enumerate}
\end{lemma}
In particular, we notice the following result from Lemma~\ref{lem:a1}:
\begin{equation}
\mathcal{S}_{D,1}[\psi](x)=-\frac{\mathrm{i}}{4\pi}\int_{\partial D}\psi(y)d\sigma(y)\quad\mbox{and}\quad \mathcal{K}_{D,1}=0.
\end{equation}

Let the boundary integral operators $\mathcal{K}_{D_\delta,j}$ and $\mathcal{S}_{D_\delta,j}$ be defined by (\ref{sing01}) and (\ref{doub01}) respectively. We next derive further asymptotic expansions of $\mathcal{K}_{D_\delta,j}$, $\mathcal{S}_{D_\delta,j}$, the single layer potential operator $\mathcal{S}_{D_\delta}$ and the static Neumann-Poincar\'{e} operator $\mathcal{K}_{D_\delta}^*$ with respect to the size scale $\delta$.

\begin{lemma}\label{lekstar01}
Let $\psi\in\mathcal {H}^*(\partial D_\delta)$ and $\tilde{\psi}(\tilde{x})=\psi(x)$ for $x\in\partial D_\delta$ and $\tilde{x}\in\partial D$, and let $\iota_{1,\delta^{t}}(\tilde{x})$, $t=\frac{1}{2}$ or $1$, be the region defined by
\begin{align}\label{}
\iota_{1,\delta^{t}}(\tilde{x}):=\{\tilde{y}\mid|z_{\tilde{x}}-z_{\tilde{y}}|<\delta^t,\tilde{y}\in\partial D\}.
\end{align}
Then we have the following asymptotic results:
\begin{enumerate}
\item[(i)]
\begin{equation}\label{8.11}
\mathcal{K}_{D_\delta}^*[\psi](x)=\mathcal{K}_{\delta,c}^*[\tilde{\psi}](\tilde{x})
+\delta\mathcal{K}_{S^f\setminus\overline{\iota_{1,\delta^{1/2}}}}^*[\tilde{\psi}](\tilde{x})+o(\delta),\ \ \ c=\{a,b\},
\end{equation}
where
\begin{align}\label{}
&\mathcal{K}_{\delta,c}^*[\tilde{\psi}](\tilde{x}):=\frac{1}{4\pi}\int_{S^c\cap\overline{\iota_{1,\delta}(\tilde{x})}}\frac{( \tilde{x}-\tilde{y}+(\delta^{-1}-1)(z_{\tilde{x}}-z_{\tilde{y}}),\nu_x)}{|\tilde{x}-\tilde{y}+(\delta^{-1}-1)(z_{\tilde{x}}-z_{\tilde{y}})|^3}
\tilde{\psi}(\tilde{y})d\sigma(\tilde{y}),\\
&\mathcal{K}_{S^f\setminus\overline{\iota_{1,\delta^{1/2}}}}^*[\tilde{\psi}](\tilde{x}):=\frac{1}{4\pi}\int_{S^f\setminus\overline{\iota_{1,\delta^{1/2}}(\tilde{x})}}\frac{( z_{\tilde{x}}-z_{\tilde{y}},\nu_x)}{|z_{\tilde{x}}-z_{\tilde{y}}|^3}
\tilde{\psi}(\tilde{y})d\sigma(\tilde{y}).
\end{align}

\item[(ii)]
\begin{equation}\label{8.14}
\mathcal{K}_{D_\delta,2}[\psi](x)=\delta^2 \mathcal{K}_{\delta,c}^{(2)}[\tilde{\psi}](\tilde{x})
+\delta \mathcal{K}_{S^f\setminus\overline{\iota_{1,\delta^{1/2}}}}^{(2)}[\tilde{\psi}](\tilde{x})+o(\delta),\ \ \ c=\{a,b\},
\end{equation}
where
\begin{align}\label{}
&\mathcal{K}_{\delta,c}^{(2)}[\tilde{\psi}](\tilde{x}):=\frac{1}{8\pi}\int_{S^c\cap\overline{\iota_{1,\delta}(\tilde{x})}}\frac{( \tilde{x}-\tilde{y}+(\delta^{-1}-1)(z_{\tilde{x}}-z_{\tilde{y}}),\nu_x)}{|\tilde{x}-\tilde{y}+(\delta^{-1}-1)(z_{\tilde{x}}-z_{\tilde{y}})|}
\tilde{\psi}(\tilde{y})d\sigma(\tilde{y}),\\
&\mathcal{K}_{S^f\setminus\overline{\iota_{1,\delta^{1/2}}}}^{(2)}[\tilde{\psi}](\tilde{x}):=\frac{1}{8\pi}\int_{S^f\setminus\overline{\iota_{1,\delta^{1/2}}(\tilde{x})}}\frac{( z_{\tilde{x}}-z_{\tilde{y}},\nu_x)}{|z_{\tilde{x}}-z_{\tilde{y}}|}
\tilde{\psi}(\tilde{y})d\sigma(\tilde{y}).
\end{align}

\item[(iii)]
\begin{equation}\label{}
(\mathcal{K}_{D_\delta,j}[\psi](x)=\delta^j\mathcal{K}_{S^c}^{(j)}[\tilde{\psi}](\tilde{x})
+\delta^{j-1}\mathcal{K}_{S^f}^{(j)}[\tilde{\psi}](\tilde{x}),\ \ \ j\geq3,\ c=\{a,b\},
\end{equation}
where
\begin{align}\label{}
&\mathcal{K}_{S^c}^{(j)}[\tilde{\psi}](\tilde{x}):=-\frac{\mathrm{i}^j(j-1)}{4\pi j!}\int_{S^c}\frac{( \tilde{x}-\tilde{y}+(\delta^{-1}-1)(z_{\tilde{x}}-z_{\tilde{y}}),\nu_x)}{|\tilde{x}-\tilde{y}+(\delta^{-1}-1)(z_{\tilde{x}}-z_{\tilde{y}})|^{3-j}}
\tilde{\psi}(\tilde{y})d\sigma(\tilde{y}),\\
&\mathcal{K}_{S^f}^{(j)}[\tilde{\psi}](\tilde{x}):=-\frac{\mathrm{i}^j(j-1)}{4\pi j!}\int_{S^f}\frac{( z_{\tilde{x}}-z_{\tilde{y}},\nu_x)}{|z_{\tilde{x}}-z_{\tilde{y}}|^{3-j}}
\tilde{\psi}(\tilde{y})d\sigma(\tilde{y}).
\end{align}
\end{enumerate}
\end{lemma}

\begin{proof}
The proofs of (i) and (ii) can follow from a similar argument to that of Lemma 4.2 in \cite{DLU7}, so we only give the proof of (iii) in what follows. From the definition of $\mathcal{K}_{D_\delta,j}$, we see that, for $j\geq 3$, the operator $\mathcal{K}_{D_\delta,j}$ has no singularity. Then by the
transformation formula (\ref{2.4}), we have
\begin{align*}\label{}
\mathcal{K}_{D_\delta,j}[\psi](x)&=-\frac{\mathrm{i}^j(j-1)}{4\pi j!}\int_{\partial D_\delta}|x-y|^{j-3}( x-y,\nu(x))\psi(y)d\sigma(y)\\
&=-\frac{\mathrm{i}^j(j-1)}{4\pi j!}\int_{S_\delta^c\cup S_\delta^f}|x-y|^{j-3}( x-y,\nu(x))\psi(y)d\sigma(y)\\
&=-\frac{\mathrm{i}^j(j-1)}{4\pi j!}\int_{S^c}\delta^j\frac{( \tilde{x}-\tilde{y}+(\delta^{-1}-1)(z_{\tilde{x}}-z_{\tilde{y}}),\nu_x)}{|\tilde{x}-\tilde{y}+(\delta^{-1}-1)(z_{\tilde{x}}-z_{\tilde{y}})|^{3-j}}
\tilde{\psi}(\tilde{y})d\sigma(\tilde{y})\\
&\ \ \ \ \ \ \ \ \ -\frac{\mathrm{i}^j(j-1)}{4\pi j!}\int_{S^f}\delta^{j-1}\frac{( \tilde{x}-\tilde{y}+(\delta^{-1}-1)(z_{\tilde{x}}-z_{\tilde{y}}),\nu_x)}{|\tilde{x}-\tilde{y}+(\delta^{-1}-1)(z_{\tilde{x}}-z_{\tilde{y}})|^{3-j}}
\tilde{\psi}(\tilde{y})d\sigma(\tilde{y})\\
&=\delta^j\cdot\mathcal{K}_{S^c}^{(j)}[\tilde{\psi}](\tilde{x})
+\delta^{j-1}\cdot\mathcal{K}_{S^f}^{(j)}[\tilde{\psi}](\tilde{x}),
\end{align*}
which yields the assertion in (iii).

The proof is complete.
\end{proof}

\begin{remark}\label{re:0101}
Since the transformation $A$ defined by (\ref{2.4}) is a diffeomorphism from $\partial D_\delta$ onto $\partial D$, we can rewrite  (\ref{8.11}) as
\begin{equation}
\mathcal{K}_{D_\delta}^*[\psi](x)=\mathcal{K}_{\delta,c}^*[A^{-1}\diamond\psi](A(x))
+\delta\mathcal{K}_{S^f\setminus\overline{\iota_{1,\delta^{1/2}}}}^*[A^{-1}\diamond\psi](A(x))+o(\delta),\ \ \ c=\{a,b\},
\end{equation}
where $(A^{-1}\diamond\psi)(\tilde{x}):=\psi(A^{-1}(\tilde{x}))=\tilde{\psi}(\tilde{x})$, known as the pull-back transformation. Thus, the integral operators $\mathcal{K}_{\delta,c}^*$ and $\mathcal{K}_{S^f\setminus\overline{\iota_{1,\delta^{1/2}}}}^*$ in the RHS of (\ref{8.11}) should be seen as the operators from $\mathcal {H}^*(\partial D_\delta)$ to $\mathcal {H}^*(\partial D_\delta)$. In what follows, we always hold this view in our asymptotic analysis and spectral expansions.
\end{remark}

We next present the expansion formulas for the operators $\mathcal{S}_{D_\delta,j}$ and $\mathcal{S}_{D_\delta}$ with respect to the size parameter $\delta\ll 1$.

\begin{lemma}\label{le:app0101}
Let $\psi\in\mathcal {H}^*(\partial D_\delta)$ and $\tilde{\psi}(\tilde{x})=\psi(x)$ for $x\in\partial D_\delta$ and $\tilde{x}\in\partial D$. Then the following asymptotic results hold:

\begin{enumerate}
\item[(i)]
\begin{equation}\label{8.20}
\mathcal{S}_{D_\delta}[\psi](x)=\delta \mathcal{S}_{\delta,c}[\tilde{\psi}](\tilde{x})
+\delta\mathcal{S}_{S^f\setminus\overline{\iota_{1,\delta^{1/2}}}}[\tilde{\psi}](\tilde{x})+o(\delta),\ \ \ c=\{a,b\},
\end{equation}
where
\begin{align}\label{}
\mathcal{S}_{\delta,c}[\tilde{\psi}](\tilde{x}):=&-\frac{1}{4\pi}\int_{S^c\cap\overline{\iota_{1,\delta}(\tilde{x})}}\frac{1}{|\tilde{x}-\tilde{y}+(\delta^{-1}-1)(z_{\tilde{x}}-z_{\tilde{y}})|}
\tilde{\psi}(\tilde{y})d\sigma(\tilde{y}),\\
\mathcal{S}_{S^f\setminus\overline{\iota_{1,\delta^{1/2}}}}[\tilde{\psi}](\tilde{x}):=&-\frac{1}{4\pi}\int_{S^f\setminus\overline{\iota_{1,\delta^{1/2}}(\tilde{x})}}\frac{ 1}{|z_{\tilde{x}}-z_{\tilde{y}}|}
\tilde{\psi}(\tilde{y})d\sigma(\tilde{y}).
\end{align}

\item[(ii)]
\begin{equation}\label{8.23}
\mathcal{S}_{D_\delta,j}[\psi](x)=\delta^{j+1}\mathcal{S}_{S^c}^{(j)}[\tilde{\psi}](\tilde{x})
+\delta^{j}\mathcal{S}_{S^f}^{(j)}[\tilde{\psi}](\tilde{x}),\ \ \ j\geq1,\ c=\{a,b\},
\end{equation}
where
\begin{align}\label{sing02}
&\mathcal{S}_{S^c}^{(j)}[\tilde{\psi}](\tilde{x}):=-\frac{\mathrm{i}^j}{4\pi j!}\int_{S^c}|\tilde{x}-\tilde{y}+(\delta^{-1}-1)(z_{\tilde{x}}-z_{\tilde{y}})|^{j-1}
\tilde{\psi}(\tilde{y})d\sigma(\tilde{y}),\\
&\mathcal{S}_{S^f}^{(j)}[\tilde{\psi}](\tilde{x}):=-\frac{\mathrm{i}^j}{4\pi j!}\int_{S^f}|\tilde{x}-\tilde{y}+(\delta^{-1}-1)(z_{\tilde{x}}-z_{\tilde{y}})|^{j-1}
\tilde{\psi}(\tilde{y})d\sigma(\tilde{y}).\label{sing03}
\end{align}
\end{enumerate}
\end{lemma}
From (\ref{sing02}) and (\ref{sing03}), it is directly verified that
\begin{equation}\label{}
\mathcal{S}_{S^c}^{(1)}[\tilde{\psi}](\tilde{x}):=-\frac{\mathrm{i}^j}{4\pi j!}\int_{S^c}
\tilde{\psi}(\tilde{y})d\sigma(\tilde{y}),\ \ \ \ \
\mathcal{S}_{S^f}^{(1)}[\tilde{\psi}](\tilde{x}):=-\frac{\mathrm{i}^j}{4\pi j!}\int_{S^f}
\tilde{\psi}(\tilde{y})d\sigma(\tilde{y}).
\end{equation}

\subsection{Asymptotics of the scattering field in the quasi-static regime}

In order to analyze the plasmon resonance for the scattering system \eqref{2.8}, we first establish the asymptotic expansion formula of the scattered field with respect to the angular frequency $\omega\ll 1$. To that end, we resent the following lemma.

\begin{lemma}[\cite{ACKLM5,AMRZ14}] \label{eigen} Let $D_\delta$ be defined in Section 2. Then
\begin{enumerate}
\item[(i)] $\mathcal K_{D_\delta}^*$ is a compact self-adjoint
operator in the Hilbert space $\mathcal {H}^*(\partial D_\delta)$.

\item[(ii)] Let $\left\{\lambda_{j,\delta};\varphi_{j,\delta}\right\},\ j=0,1,2,\cdots$, be the eigenvalue
and eigenfunction pair of $\mathcal {K}_{D_\delta}^*$,
where $\lambda_0=\frac{1}{2}$. Then,
$\lambda_{j,\delta}\in(-\frac{1}{2},\frac{1}{2}]$, and
$\lambda_{j,\delta}\rightarrow0$ as $j\rightarrow\infty$.

\item[(iii)] $\mathcal {H}^*(\partial D_{\delta})=\mathcal {H}_0^*(\partial
D_{\delta})\oplus\{c\varphi_0\},\ c\in\mathbb{C}$, where
\[
\mathcal
{H}_0^*(\partial D_{\delta})=\{\phi\in\mathcal {H}^*(\partial
D_{\delta}):\int_{\partial D_{\delta}}\phi d\sigma=0\}.
\]

\item[(iv)] For any $\psi\in\ H^{-\frac{1}{2}}(\partial D_{\delta})$, it holds that
\begin{equation}\label{3.2}
\mathcal{K}_{D_\delta}^*[\psi]=\sum_{j=0}^\infty\lambda_{j,\delta}\frac{\langle\psi,\varphi_{j,\delta}\rangle_{\mathcal
{H}^*(\partial D_{\delta})}}{\langle\varphi_{j,\delta},\varphi_{j,\delta}\rangle_{\mathcal
{H}^*(\partial D_{\delta})}}\varphi_{j,\delta}.
\end{equation}
\end{enumerate}
\end{lemma}

From now on, we use $(\cdot,\cdot)$ as the standard inner product in $\mathbb{R}^3$. The inner product (\ref{3.1}) and the corresponding norm on $\partial D_\delta$ are denoted by $\langle\cdot,\cdot\rangle$ and $\|\cdot\|$ in short, respectively. $A\lesssim B$ means $A\leq CB$ for some generic positive constant $C$. $A\approx B$ means that $A\lesssim B$ and $B\lesssim A$.

Owing to (\ref{3.2}) and (\ref{2.20}), we have
\begin{equation}\label{3.3}
\mathcal
{A}_{D_\delta,0}[\psi]=\sum_{j=0}^\infty\tau_{j,\delta}\frac{\langle\psi,\varphi_{j,\delta}\rangle}
{\langle\varphi_{j,\delta},\varphi_{j,\delta}\rangle}\varphi_{j,\delta},
\end{equation}
where
\begin{equation}\label{3.4}
\tau_{j,\delta}=\frac{1}{2}\left(\frac{1}{\varepsilon_m}+\frac{1}{\varepsilon_c}\right)+\left(\frac{1}{\varepsilon_m}-\frac{1}{\varepsilon_c}\right)\lambda_{j,\delta}.
\end{equation}
In what follows, for notational convenience, we define
\begin{equation}\label{eq:n1}
a_{j,\delta}:=\langle\varphi_{j,\delta},\varphi_{j,\delta}\rangle_{\mathcal
{H}^*(\partial D_{\delta})}.
\end{equation}

We present the following lemma on the asymptotic expansion of the operator $\mathcal{A}_{D_\delta}(\omega)$ introduced in \eqref{2.18} with respect to $\omega\ll 1$, whose proof follows from straightforward computations and is omitted.
\begin{lemma}\label{lem3.2}
The operator $\mathcal{A}_{D_\delta}(\omega):\ \mathcal {H}^*(\partial
D_\delta)\rightarrow\mathcal {H}^*(\partial D_\delta)$ has the following expansion:
\begin{equation}\label{3.5}
\mathcal {A}_{D_\delta}(\omega)=\mathcal {A}_{D_\delta,0}+\omega^{2}\mathcal
{A}_{D_\delta,2}+\mathcal{O}(\omega^3)
\end{equation}
with
\begin{equation}\label{AD2}
\mathcal {A}_{D_\delta,2}=(\mu_m-\mu_c)\mathcal
{K}_{D_\delta,2}+\frac{\varepsilon_m\mu_m-\varepsilon_c\mu_c}{\varepsilon_c}\left(\frac{1}{2}\mathcal
{I}-\mathcal {K}_{D_\delta}^*\right)\mathcal{S}_{D_\delta}^{-1}\mathcal{S}_{D_\delta,2},
\end{equation}
where $\mathcal{K}_{D_\delta,2}$ and $\mathcal{S}_{D_\delta,2}$ are defined in Section \ref{sec2}.
\end{lemma}

For our subsequent use, we introduce the so-called index set of plasmon resonance.

\begin{definition}\label{de2}
Let $\tau_{j,\delta}$ be introduced in \eqref{3.4}. We say that $J\subseteq \mathbb{N}$ is an index set of
resonance if $\tau_{j,\delta}$ is close to zero when $j\in J$ and is
bounded from below when $j\in J^c:=\mathbb{N}\backslash J$. More precisely, we choose a
threshold number $\eta_0>0$ independent of $\omega$ and $\delta$ such that
$|\tau_{j,\delta}|\geq\eta_0>0$, for $j\in J^c$.
\end{definition}

Next, we impose the following two mild conditions throughout our study:
\begin{enumerate}
\item[(C1)] Each eigenvalue $\lambda_{j,\delta}$ for $j\in J$ is a simple
eigenvalue of the operator $\mathcal {K}_{D_\delta}^*$.
\item[(C2)] Suppose that $\varepsilon_c+\varepsilon_m\neq0$.
\end{enumerate}
It is noted that in \eqref{2.6} and \eqref{2.7}, we assume that $\Im \varepsilon_c>0$, whereas $\varepsilon_m$ is a real constant. Hence, condition (C2) is easily fulfilled.
We would like to point out that by condition (C2), one can deduce that the index set $J$ is finite. Noting also that for $j=0$, $\lambda_{0,\delta}=\frac{1}{2}$, we see
that $\tau_{0,\delta}=\frac{1}{\mu_m}\sim 1$. Thus, throughout this paper, we exclude $0$ from the index
set $J$.

\begin{lemma}\label{freasy01}
In the quasi-static regime and under conditions (C1) and (C2), the scattering field $u^s$ to \eqref{2.8} has
the following representation
\begin{equation}\label{3.6}
u^s=\mathcal {S}_{D_\delta}^{k_m}[\psi],
\end{equation}
where
\begin{equation}\label{3.7}
\psi=\sum_{j\in J}\frac{\mathrm{i}\omega\sqrt{\mu_m\varepsilon_m}(1/\varepsilon_c-1/\varepsilon_m)a_{j,\delta}^{-1}\langle d\cdot\nu,\varphi_{j,\delta}\rangle\varphi_{j,\delta}+\mathcal{O}(\omega^2)}{\tau_{j,\delta}+\mathcal{O}(\omega^2)}+\mathcal{O}(\omega).
\end{equation}
\end{lemma}
\begin{proof}
Note that the incident wave $u^i=e^{\mathrm{i}k_m d\cdot x}$ admits the following asymptotic expansion:
\beq
u^i=1+\mathrm{i}k_m d\cdot x +\mathcal{O}(\omega^2).
\eeq
By combing \eqnref{2.17}, \eqnref{2.18}, \eqnref{2.20}, together with the spectral expansion \eqnref{3.3}, and the fact that
\beq
\left(\frac{\mathcal {I}}{2}-\mathcal
{K}_{D_\delta}^*\right)\mathcal {S}_{D_\delta}^{-1}[u^i]=\mathcal {S}_{D_\delta}^{-1}\left(\frac{\mathcal {I}}{2}-\mathcal
{K}_{D_\delta}\right)[u^i]=\mathrm{i}k_m {S}_{D_\delta}^{-1}\left(\frac{\mathcal {I}}{2}-\mathcal
{K}_{D_\delta}\right)[x\cdot d]+\mathcal{O}(\omega^2)
\eeq
 one can obtain \eqnref{3.7}.
\end{proof}

\subsection{Asymptotics of the scattering field for $\delta\ll 1$}

\eqref{3.6}-\eqref{3.7} give the asymptotic expansion of the scattered wave $u^s$ with respect to $\omega\ll 1$. It is noted that the anisotropic size parameter $\delta$ is also asymptotically small.
In this part, we shall derive further asymptotic expansions with respect to the size parameter $\delta$. We stress that this part is essential in our analysis, since the nanorod is anisotropic with respect to its dimensional sizes. First, by using the expansions of the layer potential operators $\mathcal{S}_{D_\delta}$ and $\mathcal{K}_{D_\delta}^*$ with respect to $\delta$, we have the following lemma.

\begin{lemma}\label{3.8}
The operator $\mathcal {A}_{D_\delta}(\omega):\ \mathcal {H}^*(\partial
D_\delta)\rightarrow\mathcal {H}^*(\partial D_\delta)$ has the expansion
formula as follows
\begin{equation}\label{3.8}
\mathcal {A}_{D_\delta}(\omega)=\mathcal {A}_{D_\delta,0}+\omega^{2}\delta\widehat{\mathcal
{A}}_{D_\delta,2}+o(\omega^2\delta)+\mathcal{O}(\omega^3),
\end{equation}
where
\[
\widehat{\mathcal {A}}_{D_\delta,2}=(\mu_m-\mu_c)\mathcal{K}_{S^f\setminus\overline{\iota_{1,\delta^{1/2}}}}^{(2)}
+\frac{\varepsilon_m\mu_m-\varepsilon_c\mu_c}{\varepsilon_c}\left(\frac{1}{2}\mathcal
{I}-\mathcal{K}_{\delta,c}^*\right)\left(\mathcal{S}_{\delta,c}+\mathcal{S}_{S^f\setminus\overline{\iota_{1,\delta^{1/2}}}}
+o(1)\right)^{-1}\mathcal{S}_{S^f}^{(2)}
\]
with $\mathcal{K}_{\delta,c}^*$, $\mathcal{K}_{S^f\setminus\overline{\iota_{1,\delta^{1/2}}}}^{(2)}$, $\mathcal{S}_{\delta,c}$, $\mathcal{S}_{S^f\setminus\overline{\iota_{1,\delta^{1/2}}}}$ and $\mathcal{S}_{S^f}^{(2)}$ defined in Section \ref{sec2}.
\end{lemma}

\begin{proof}
Since $\mathcal{S}_{D_\delta}$ is invertible, by
substituting the expansion formulas (\ref{8.11}), (\ref{8.14}), (\ref{8.20}) and (\ref{8.23}) into (\ref{3.5}), one can derive (\ref{3.8}) by direct calculations.
\end{proof}

For the subsequent use, we define the following regions associated with $D_\delta$ and $D$:
\begin{align}\label{}
\iota_{\delta}(P_0):=\{y\mid|P_0-z_{y}|<\delta,y\in\partial D_\delta\},\\
\iota_{\delta}(Q_0):=\{y\mid|Q_0-z_{y}|<\delta,y\in\partial D_\delta\},\\
\iota_{1,\delta}(P_0):=\{\tilde{y}\mid|P_0-z_{\tilde{y}}|<\delta,\tilde{y}\in\partial D\},\\
\iota_{1,\delta}(Q_0):=\{\tilde{y}\mid|Q_0-z_{\tilde{y}}|<\delta,\tilde{y}\in\partial D\}.
\end{align}
In what follows, we use $\|\cdot\|_\Gamma$ to represent the norm in $\mathcal {H}^*(\Gamma)$ (cf. (\ref{3.1})) for any surface $\Gamma$. The following lemma is of critical importance for our subsequent analysis.
\begin{lemma}\label{le:3.6}
Let $\mathcal{K}_{\delta,c}^*$ and $\mathcal{K}_{S^f\setminus\overline{\iota_{1,\delta^{1/2}}}}^*$ be defined in Lemma \ref{lekstar01}. Then the following results hold.

\begin{enumerate}

\item[(i)] If $x\in S_\delta^a$ or ($\tilde{x}\in S^a$), for any $\tilde{\psi}\in\mathcal {H}^*(\partial
D)$, we have
\begin{align}\label{}
&\mathcal{K}_{\delta,a}^*[\tilde{\psi}](\tilde{x})=\mathcal{K}_{S^a}^*[\tilde{\psi}](\tilde{x}):=\frac{1}{4\pi}\int_{S^a}\frac{( \tilde{x}-\tilde{y},\nu_x)}{|\tilde{x}-\tilde{y}|^3}\tilde{\psi}(\tilde{y})d\sigma(\tilde{y}),\\
&\mathcal{K}_{S^f\setminus\overline{\iota_{1,\delta^{1/2}}}}^*[\tilde{\psi}](\tilde{x})=\mathcal{K}_{S^f,P_0}^*[\tilde{\psi}](\tilde{x})
+o\left(\delta^{\frac{1}{2}}\|\tilde\psi\|_{\partial D}\right),
\end{align}
where
\begin{align}\label{}
\mathcal{K}_{S^f,P_0}^*[\tilde{\psi}](\tilde{x})=\frac{1}{4\pi}\int_{S^f}\frac{( P_0-z_{\tilde{y}},\nu_{P_0})}{|P_0-z_{\tilde{y}}|^3}\tilde{\psi}(\tilde{y})d\sigma(\tilde{y}).
\end{align}

\item[(ii)] If $x\in S_\delta^b$ or ($\tilde{x}\in S^b$), for any $\tilde{\psi}\in\mathcal {H}^*(\partial
D)$, we have
\begin{align}\label{}
&\mathcal{K}_{\delta,b}^*[\tilde{\psi}](\tilde{x})=\mathcal{K}_{S^b}^*[\tilde{\psi}](\tilde{x}):=\frac{1}{4\pi}\int_{S^b}\frac{( \tilde{x}-\tilde{y},\nu_x)}{|\tilde{x}-\tilde{y}|^3}\tilde{\psi}(\tilde{y})d\sigma(\tilde{y}),\\
&\mathcal{K}_{S^f\setminus\overline{\iota_{1,\delta^{1/2}}}}^*[\tilde{\psi}](\tilde{x})=\mathcal{K}_{S^f,Q_0}^*[\tilde{\psi}](\tilde{x})
+o\left(\delta^{\frac{1}{2}}\|\tilde\psi\|_{\partial D}\right),
\end{align}
where
\begin{align}\label{}
\mathcal{K}_{S^f,Q_0}^*[\tilde{\psi}](\tilde{x})=\frac{1}{4\pi}\int_{S^f}\frac{( Q_0-z_{\tilde{y}},\nu_{Q_0})}{|Q_0-z_{\tilde{y}}|^3}\tilde{\psi}(\tilde{y})d\sigma(\tilde{y}).
\end{align}

\item[(iii)] If $x\in S_\delta^f\cap\overline{\iota_\delta(P_0)}$ or ($\tilde{x}\in S^f\cap\overline{\iota_{1,\delta}(P_0)}$), for any $\tilde{\psi}\in\mathcal {H}^*(\partial
D)$, we have
\begin{align}\label{}
&\mathcal{K}_{\delta,a}^*[\tilde{\psi}](\tilde{x})=\mathcal{K}_{S^a,P_0}^*[\tilde{\psi}](\tilde{x})+\mathcal{O}\left(\delta\|\tilde\psi\|_{\partial D}\right),\\
&\mathcal{K}_{S^f\setminus\overline{\iota_{1,\delta^{1/2}}}}^*[\tilde{\psi}](\tilde{x})=\mathcal{K}_{S^f}^*[\tilde{\psi}](\tilde{x})
+\mathcal{O}\left(\delta^{\frac{1}{2}}\|\tilde\psi\|_{\partial D}\right),
\end{align}
where
\begin{align}\label{eq:defKSAP01}
&\mathcal{K}_{S^a,P_0}^*[\tilde{\psi}](\tilde{x})=\frac{1}{4\pi}\int_{S^a}\frac{( \tilde{x}_1-\tilde{y},\nu_{\tilde{x}_1})}{|\tilde{x}_1-\tilde{y}|^3}\tilde{\psi}(\tilde{y})d\sigma(\tilde{y}),\ \ \ \tilde{x}_1\in \overline{S^a}\cap\overline{S^f},\\
&\mathcal{K}_{S^f}^*[\tilde{\psi}](\tilde{x})=\frac{1}{4\pi}\int_{S^f}\frac{( z_{\tilde{x}}-z_{\tilde{y}},\nu_{x})}{|z_{\tilde{x}}-z_{\tilde{y}}|^3}\tilde{\psi}(\tilde{y})d\sigma(\tilde{y}).
\end{align}

\item[(iv)] If $x\in S_\delta^f\cap\overline{\iota_\delta(Q_0)}$ or ($\tilde{x}\in S^f\cap\overline{\iota_{1,\delta}(Q_0)}$), for any $\tilde{\psi}\in\mathcal {H}^*(\partial
D)$, we have
\begin{align}\label{}
&\mathcal{K}_{\delta,b}^*[\tilde{\psi}](\tilde{x})=\mathcal{K}_{S^b,Q_0}^*[\tilde{\psi}](\tilde{x})+\mathcal{O}\left(\delta\|\tilde\psi\|_{\partial D}\right),\\
&\mathcal{K}_{S^f\setminus\overline{\iota_{1,\delta^{1/2}}}}^*[\tilde{\psi}](\tilde{x})=\mathcal{K}_{S^f}^*[\tilde{\psi}](\tilde{x})
+\mathcal{O}\left(\delta^{\frac{1}{2}}\|\tilde\psi\|_{\partial D}\right),
\end{align}
where
\begin{align}\label{eq:defKSBQ01}
\mathcal{K}_{S^b,Q_0}^*[\tilde{\psi}](\tilde{x})=\frac{1}{4\pi}\int_{S^b}\frac{( \tilde{x}_1-\tilde{y},\nu_{\tilde{x}_1})}{|\tilde{x}_1-\tilde{y}|^3}\tilde{\psi}(\tilde{y})d\sigma(\tilde{y}),\ \ \ \tilde{x}_1\in \overline{S^b}\cap\overline{S^f}.
\end{align}

\item[(v)] If $x\in S_\delta^f\setminus\overline{\iota_\delta(P_0)\cup\iota_\delta(Q_0)}$ or ($\tilde{x}\in S^f\setminus\overline{\iota_{1,\delta}(P_0)\cup\iota_{1,\delta}(Q_0)}$), for any $\tilde{\psi}\in\mathcal {H}^*(\partial
D)$, we have
\begin{align}\label{}
&\mathcal{K}_{\delta,c}^*[\tilde{\psi}](\tilde{x})=0,\\
&\mathcal{K}_{S^f\setminus\overline{\iota_{1,\delta^{1/2}}}}^*[\tilde{\psi}](\tilde{x})=\mathcal{K}_{S^f}^*[\tilde{\psi}](\tilde{x})
+\mathcal{O}\left(\delta^{\frac{1}{2}}\|\tilde\psi\|_{\partial D}\right).
\end{align}
\end{enumerate}

\end{lemma}

\begin{proof}
Since the proofs of (i)-(v) are similar, we only prove (i) in what follows.
By the definition of $\mathcal{K}_{\delta,a}^*$ and noting that $S^a\cap\iota_{1,\delta}(\tilde{x})=S^a$ if $\tilde{x}\in S^a$, one has
\begin{align*}\label{}
\mathcal{K}_{\delta,a}^*[\tilde{\psi}](\tilde{x})&=\frac{1}{4\pi}\int_{S^a\cap\overline{\iota_{1,\delta}(\tilde{x})}}\frac{( \tilde{x}-\tilde{y}+(\delta^{-1}-1)(z_{\tilde{x}}-z_{\tilde{y}}),\nu_x)}{|\tilde{x}-\tilde{y}+(\delta^{-1}-1)(z_{\tilde{x}}-z_{\tilde{y}})|^3}
\tilde{\psi}(\tilde{y})d\sigma(\tilde{y})\\
&=\frac{1}{4\pi}\int_{S^a}\frac{(\tilde{x}-\tilde{y},\nu_x)}{|\tilde{x}-\tilde{y}|^3}\tilde{\psi}(\tilde{y})d\sigma(\tilde{y}).
\end{align*}
Moreover by direct asymptotic analysis, one has
\begin{align*}\label{}
\mathcal{K}_{S^f\setminus\overline{\iota_{1,\delta^{1/2}}}}^*[\tilde{\psi}](\tilde{x})
&=\frac{1}{4\pi}\int_{S^f\setminus\overline{\iota_{1,\delta^{1/2}}(\tilde{x})}}\frac{( z_{\tilde{x}}-z_{\tilde{y}},\nu_{x})}{|z_{\tilde{x}}-z_{\tilde{y}}|^3}\tilde{\psi}(\tilde{y})d\sigma(\tilde{y})\\
&=\frac{1}{4\pi}\int_{S^f}\frac{( P_0-z_{\tilde{y}},\nu_{P_0})}{|P_0-z_{\tilde{y}}|^3}\tilde{\psi}(\tilde{y})d\sigma(\tilde{y})-\frac{1}{4\pi}
\int_{\overline{\iota_{1,\delta^{1/2}}(P_0)}}\frac{( z_{\tilde{x}}-z_{\tilde{y}},\nu_{x})}{|z_{\tilde{x}}-z_{\tilde{y}}|^3}\tilde{\psi}(\tilde{y})d\sigma(\tilde{y})\\
&=\frac{1}{4\pi}\int_{S^f}\frac{( P_0-z_{\tilde{y}},\nu_{P_0})}{|P_0-z_{\tilde{y}}|^3}\tilde{\psi}(\tilde{y})d\sigma(\tilde{y})+o\left(\delta^{\frac{1}{2}}\|\tilde\psi\|_{\partial D}\right).
\end{align*}

The proof is complete.
\end{proof}

In a similar manner, one can derive the expansion formulas of $\mathcal{S}_{\delta,c}$ and
$\mathcal{S}_{S^f\setminus\overline{\iota_{1,\delta^{1/2}}}}$ as follows.

\begin{lemma}\label{le:sigl01}
Let $\mathcal{S}_{\delta,c}$, $\mathcal{S}_{S^f\setminus\overline{\iota_{1,\delta^{1/2}}}}$ be defined in Lemma \ref{le:app0101}. Then the following results hold.

(i) If $x\in S_\delta^a$ (or equivalently $\tilde{x}\in S^a$), for any $\tilde{\psi}\in\mathcal {H}^*(\partial
D)$, we have
\begin{align}\label{}
&\mathcal{S}_{\delta,a}[\tilde{\psi}](\tilde{x})=\mathcal{S}_{S^a}[\tilde{\psi}](\tilde{x}):=-\frac{1}{4\pi}\int_{S^a}\frac{1}{|\tilde{x}-\tilde{y}|}\tilde{\psi}(\tilde{y})d\sigma(\tilde{y}),\\
&\mathcal{S}_{S^f\setminus\overline{\iota_{1,\delta^{1/2}}}}[\tilde{\psi}](\tilde{x})=\mathcal{S}_{S^f,P_0}[\tilde{\psi}](\tilde{x})
+o\left(\delta^{\frac{1}{2}}\|\tilde\psi\|_{\partial D}\right),
\end{align}
where
\begin{align}\label{}
\mathcal{S}_{S^f,P_0}[\tilde{\psi}](\tilde{x})=-\frac{1}{4\pi}\int_{S^f}\frac{1}{|P_0-z_{\tilde{y}}|}\tilde{\psi}(\tilde{y})d\sigma(\tilde{y}).
\end{align}

(ii) If $x\in S_\delta^b$ (or equivalently $\tilde{x}\in S^b$), for any $\tilde{\psi}\in\mathcal {H}^*(\partial
D)$, we have
\begin{align}\label{}
&\mathcal{S}_{\delta,b}[\tilde{\psi}](\tilde{x})=\mathcal{S}_{S^b}[\tilde{\psi}](\tilde{x}):
=-\frac{1}{4\pi}\int_{S^b}\frac{1}{|\tilde{x}-\tilde{y}|}\tilde{\psi}(\tilde{y})d\sigma(\tilde{y}),\\
&\mathcal{S}_{S^f\setminus\overline{\iota_{1,\delta^{1/2}}}}[\tilde{\psi}](\tilde{x})=\mathcal{S}_{S^f,Q_0}[\tilde{\psi}](\tilde{x})
+o\left(\delta^{\frac{1}{2}}\|\tilde\psi\|_{\partial D}\right),
\end{align}
where
\begin{align}\label{}
\mathcal{S}_{S^f,Q_0}[\tilde{\psi}](\tilde{x})=-\frac{1}{4\pi}\int_{S^f}\frac{1}{|Q_0-z_{\tilde{y}}|}\tilde{\psi}(\tilde{y})d\sigma(\tilde{y}).
\end{align}

(iii) If $x\in S_\delta^f\cap\overline{\iota_\delta(P_0)}$ (or equivalently $\tilde{x}\in S^f\cap\overline{\iota_{1,\delta}(P_0)}$), for any $\tilde{\psi}\in\mathcal {H}^*(\partial
D)$, we have
\begin{align}
&\mathcal{S}_{\delta,a}[\tilde{\psi}](\tilde{x})=\mathcal{S}_{S^a,P_0}[\tilde{\psi}](\tilde{x})+\mathcal{O}\left(\delta\|\tilde\psi\|_{\partial D}\right),\label{eq:nn1}\\
&\mathcal{S}_{S^f\setminus\overline{\iota_{1,\delta^{1/2}}}}[\tilde{\psi}](\tilde{x})=\mathcal{S}_{S^f}[\tilde{\psi}](\tilde{x})
+\mathcal{O}\left(\delta^{\frac{1}{2}}\|\tilde\psi\|_{\partial D}\right),\label{eq:nn2}
\end{align}
where
\begin{align}\label{}
&\mathcal{S}_{S^a,P_0}[\tilde{\psi}](\tilde{x})=-\frac{1}{4\pi}\int_{S^a}\frac{1}{|\tilde{x}_1-\tilde{y}|}\tilde{\psi}(\tilde{y})d\sigma(\tilde{y}),\ \ \ \tilde{x}_1\in \overline{S^a}\cap\overline{S^f},\\
&\mathcal{S}_{S^f}[\tilde{\psi}](\tilde{x})=-\frac{1}{4\pi}\int_{S^f}\frac{1}{|z_{\tilde{x}}-z_{\tilde{y}}|}\tilde{\psi}(\tilde{y})d\sigma(\tilde{y}).
\end{align}
Furthermore, the higher order terms in \eqref{eq:nn1} and \eqref{eq:nn2} are boundary integrals on subsets of $S^f$.

(iv) If $x\in S_\delta^f\cap\overline{\iota_\delta(Q_0)}$ (or equivalently $\tilde{x}\in S^f\cap\overline{\iota_{1,\delta}(Q_0)}$), for any $\tilde{\psi}\in\mathcal {H}^*(\partial
D)$, we have
\begin{align}\label{}
&\mathcal{S}_{\delta,b}[\tilde{\psi}](\tilde{x})=\mathcal{S}_{S^b,Q_0}[\tilde{\psi}](\tilde{x})+\mathcal{O}\left(\delta\|\tilde\psi\|_{\partial D}\right),\label{eq:nn3}\\
&\mathcal{S}_{S^f\setminus\overline{\iota_{1,\delta^{1/2}}}}[\tilde{\psi}](\tilde{x})=\mathcal{S}_{S^f}[\tilde{\psi}](\tilde{x})
+\mathcal{O}\left(\delta^{\frac{1}{2}}\|\tilde\psi\|_{\partial D}\right),\label{eq:nn4}
\end{align}
where
\begin{align}\label{}
\mathcal{S}_{S^b,Q_0}[\tilde{\psi}](\tilde{x})=-\frac{1}{4\pi}\int_{S^b}\frac{1}{|\tilde{x}_1-\tilde{y}|}\tilde{\psi}(\tilde{y})d\sigma(\tilde{y}),\ \ \ \tilde{x}_1\in \overline{S^b}\cap\overline{S^f}.
\end{align}
Furthermore, the higher order terms in \eqref{eq:nn3} and \eqref{eq:nn4} are boundary integrals on subsets of $S^f$.

(v) If $x\in S_\delta^f\setminus\overline{\iota_\delta(P_0)\cup\iota_\delta(Q_0)}$ (or equivalently $\tilde{x}\in S^f\setminus\overline{\iota_{1,\delta}(P_0)\cup\iota_{1,\delta}(Q_0)}$), for any $\tilde{\psi}\in\mathcal {H}^*(\partial
D)$, we have
\begin{align}\label{}
&\mathcal{S}_{\delta,c}[\tilde{\psi}](\tilde{x})=0,\\
&\mathcal{S}_{S^f\setminus\overline{\iota_{1,\delta^{1/2}}}}[\tilde{\psi}](\tilde{x})=\mathcal{S}_{S^f}[\tilde{\psi}](\tilde{x})
+\mathcal{O}\left(\delta^{\frac{1}{2}}\|\tilde\psi\|_{\partial D}\right).
\end{align}
\end{lemma}

Using the asymptotic results in Lemmas~\ref{le:3.6} and \ref{le:sigl01} and comparing with the asymptotic formula of Neumann-Poincar\'{e} operator $\mathcal {K}_{D_\delta}^*$ in (\ref{8.11}), we can obtain the following more refined one, such that the asymptotic operators at the right hand side of (\ref{8.11}) are independent of size scale $\delta$.
\begin{theorem}\label{th:NP01}
Let $\psi\in\mathcal {H}^*(\partial D_\delta)$ and $\tilde{\psi}(\tilde{x})=\psi(x)$ for $x\in\partial D_\delta$ and $\tilde{x}\in\partial D$,
Then, we have that
\begin{align}\label{Kexpan}
\mathcal{K}_{D_\delta}^*[\psi](x)=\mathcal{K}_{0}^*[\tilde{\psi}](\tilde{x})
+\delta\mathcal{K}_{1}^*[\tilde{\psi}](\tilde{x})+o(\delta\|\tilde\psi\|_{\partial D}),
\end{align}
where
\begin{equation}\label{eq:K0def01}
\begin{split}
\mathcal{K}_{0}^*[\tilde{\psi}](\tilde{x}):=&\chi\left(S^a\right)\mathcal{K}_{S^a}^*[\tilde{\psi}](\tilde{x})
+\chi\left(S^f\cap\overline{\iota_{1,\delta}(P_0)}\right)\mathcal{K}_{S^a,P_0}^*[\tilde{\psi}](\tilde{x})\\
&\ \ \ \ +\chi\left(S^f\cap\overline{\iota_{1,\delta}(Q_0)}\right)\mathcal{K}_{S^b,Q_0}^*[\tilde{\psi}](\tilde{x})
+\chi\left(S^b\right)\mathcal{K}_{S^b}^*[\tilde{\psi}](\tilde{x}),
\end{split}
\end{equation}
and
\[
\mathcal{K}_{1}^*[\tilde{\psi}](\tilde{x}):=\chi\left(S^a\right)\mathcal{K}_{S^f,P_0}^*[\tilde{\psi}](\tilde{x})
+\chi\left(S^f\right)\mathcal{K}_{S^f}^*[\tilde{\psi}](\tilde{x})+\chi\left(S^b\right)\mathcal{K}_{S^f,Q_0}^*[\tilde{\psi}](\tilde{x}).
\]
Here, $\chi$ denotes the characteristic function and the asymptotic operators $\mathcal{K}_{S^a}^*$, $\mathcal{K}_{S^b}^*$, $\mathcal{K}_{S^f}^*$, $\mathcal{K}_{S^a,P_0}^*$, $\mathcal{K}_{S^b,Q_0}^*$, $\mathcal{K}_{S^f,P_0}^*$ and $\mathcal{K}_{S^f,Q_0}^*$ are defined in
Lemma \ref{le:3.6}.
\end{theorem}
\begin{proof}
The proof follows from straightforward though a bit tedious calculations along with the use of Lemma \ref{le:3.6} and (\ref{8.11}).
\end{proof}

In what follows, we define $\widetilde{\varphi}_{j,\delta}(\tilde{x}):=\varphi_{j,\delta}(x)$ for $x\in\partial D_\delta$, and $\tilde{x}\in \partial D$.
By using Lemma \ref{eigen} and the anisotropic geometry of the nanorod $D_\delta$, one can derive the following result.
\begin{lemma}\label{le:0101}
The eigenfunction $\varphi_{j,\delta}$ in \eqnref{3.2} fulfils that
\beq\label{eq:bdint01}
\begin{split}
\int_{S^a\cup S^b}\widetilde{\varphi}_{j,\delta}=0, \quad \mbox{and} \quad \int_{S^f}\widetilde{\varphi}_{j,\delta}=0, \quad j\neq 0,
\end{split}
\eeq
\end{lemma}
\begin{proof}
Note that $\{\varphi_{j,\delta}\}_{j=0}^{\infty}$ are the eigenfunctions of $\mathcal{K}_{D_\delta}^*$ in $\mathcal{H}^*(\partial D_{\delta})$. Direct computation shows
\beq
\int_{\partial D_\delta} \varphi_{j,\delta}=-\langle \mathcal{S}_{D_\delta}^{-1}[1], \varphi_{j,\delta}\rangle=0, \quad \mbox{for} \quad j\neq 0.
\eeq
By further using the asymptotic expansion
\beq
\int_{\partial D_\delta} \varphi_{j,\delta}=\int_{S_\delta^a\cup S_\delta^b} \varphi_{j,\delta}+\int_{S_\delta^f} \varphi_{j,\delta}=\delta^2\int_{S^a\cup S^b}  \widetilde\varphi_{j,\delta}+\delta\int_{S^f}  \widetilde\varphi_{j,\delta},
\eeq
one immediately achieves \eqnref{eq:bdint01}.

The proof is complete.
\end{proof}

Let $P$ denote the eigenprojection of $\mathcal{K}_{0}^*$ on $\mathcal{H}^*(\partial D_\delta)$, i.e., $P:\mathcal{H}^*(\partial D_\delta)\rightarrow V_{\lambda_j}$, where $V_{\lambda_j}$ is the eigenspace associated with the eigenvalue $\lambda_j$ of $\mathcal{K}_{0}^*$. Furthermore, we define the reduced resolvent for $\lambda_j$ of $\mathcal{K}_{0}^*$ as (see \cite{K38})
\beq\label{eq:redureso}
\Lambda=\frac{1}{2\pi \mathrm{i}}\int_{\gamma_j}\frac{(\mathcal{K}_{0}^*-\xi)^{-1}}{\xi-\lambda_j}d\xi,
\eeq
where $\gamma_j: |\xi-\lambda_j|=r$ is a circle in the complex plane enclosing the isolated eigenvalue $\lambda_j$.

By using the perturbation theory for the spectrum of a linear operator, one can obtain the following elementary result.
\begin{lemma}\label{le:0102}
The operator $\mathcal{K}_{0}^*$ defined in \eqnref{eq:K0def01} admits the following spectral property,
\beq\label{eq:eigenasy00}
\mathcal{K}_{0}^*[\widetilde\varphi_{j}](\tilde{x})=\lambda_j \widetilde\varphi_{j}(\tilde{x})
\eeq
with $\widetilde\varphi_{j}(\tilde{x}):=\varphi_{j}(x)$ and $\lambda_j$ and $\varphi_{j}$ satisfying
\begin{equation}\label{eq:eigenasy01}
\lambda_{j,\delta}=\lambda_{j}+\delta\lambda_{j,1}+o(\delta), \quad
\varphi_{j,\delta}=\varphi_{j}-\delta\Lambda\mathcal{K}_{1}^*[\varphi_{j}]+o(\delta),
\end{equation}
where $\lambda_{j,1}$ is the eigenvalue of the operator $P(\mathcal{K}_{1}^*)P$ considered
in the eigenspace $V_{\lambda_j}$. Furthermore, one has
\beq\label{eq:eigenasy02}
\widetilde\varphi_j(\tilde{x})=0, \quad \tilde{x}\in S^f.
\eeq
\end{lemma}
\begin{proof}
Recall that $\lambda_{j,\delta}$ and $\varphi_{j,\delta}$ are the eigenvalues and eigenfunctions of $\mathcal{K}_{D_\delta}^*$ in $\mathcal{H}^*(\partial \delta)$, i.e.,
$$
\mathcal{K}_{D_\delta}^*[\varphi_{j,\delta}]=\lambda_{j,\delta}\varphi_{j,\delta}.
$$
By Theorem \ref{th:NP01} and the perturbation theory for the spectrum of linear operators (see \cite{K38}, pages 445-446, Theorem 2.6), one can find that
the eigenvalues $\lambda_{j,\delta}$ and eigenfunctions $\varphi_{j,\delta}$ of $\mathcal{K}_{D_\delta}^*$ can be approximated in \eqnref{eq:eigenasy01}, where
$\lambda_j$ and $\varphi_{j}$ satisfy \eqnref{eq:eigenasy00}.
From the definition of $\mathcal{K}_{0}^*$ in \eqnref{eq:K0def01}, one can easily see that
$$
\widetilde\varphi_j=0 \quad \mbox{in} \quad S^f\setminus\overline{\iota_{1,\delta}(P_0)\cup\iota_{1,\delta}(Q_0)}.
$$
Furthermore, if $\tilde{x}\in S^f\cap \overline{\iota_{1,\delta}(P_0)}$, then one has
$$
\mathcal{K}_{S^a,P_0}^*[\widetilde{\varphi}_{j}](\tilde{x})=\mathcal{K}_{0}^*[\widetilde{\varphi}_j](\tilde{x})=\lambda_j \widetilde\varphi_j(\tilde{x}),
$$
and by the definition of $\mathcal{K}_{S^a,P_0}^*$ in \eqnref{eq:defKSAP01}, one thus has
$$
\widetilde\varphi_j(\tilde{x})=0, \quad \tilde{x}\in S^f\cap \overline{\iota_{1,\delta}(P_0)}.
$$
Similarly, one has
$$
\widetilde\varphi_j(\tilde{x})=0, \quad \tilde{x}\in S^f\cap \overline{\iota_{1,\delta}(Q_0)}.
$$
One thus has \eqnref{eq:eigenasy02} and completes the proof.
\end{proof}

\begin{remark}
In fact, noting the definition of $\mathcal{K}_{1}^*$ and using (\ref{eq:eigenasy02}), we find $\mathcal{K}_{1}^*[\varphi_{j}]=0$, and then the expansion formula of the eigenfunction in (\ref{eq:eigenasy01}) can be rewritten as
\begin{equation}\label{eq:eigenasy03}
\widetilde\varphi_{j,\delta}=\widetilde\varphi_{j}+o(\delta).
\end{equation}
\end{remark}

By using (\ref{Kexpan}), together with Lemmas \ref{le:0101} and \ref{le:0102}, the expression of the layer potential density associated with the nanorod can be deduced as follows,
\begin{theorem}\label{th:density01}
{Under the conditions (C1) and (C2), and assume that $\psi\in\mathcal {H}^*(\partial D_\delta)$ is the solution of the operator equation (\ref{2.17}), where $\psi(x)=\widetilde{\psi}(\tilde{x})$ for $\tilde{x}\in\partial D$. Then, for every $x\in\partial D_\delta$, the solution $\psi$ can be expressed as}
$\psi(x)=\psi_1(x)+\mathcal{O}(\omega)$, where $\psi_1$ is defined by
\begin{equation}\label{3.48}
\psi_1(x)=\sum_{j\in J}\frac{\mathrm{i}\omega\sqrt{\mu_m\varepsilon_m}(1/\varepsilon_c-1/\varepsilon_m)a_{j,\delta}^{-1}\langle d\cdot\nu,\varphi_{j,\delta}\rangle\widetilde\varphi_{j}(\tilde{x})+o(\omega\delta)}
{\frac{1}{2}\left(\frac{1}{\varepsilon_m}+\frac{1}{\varepsilon_c}\right)+\left(\frac{1}{\varepsilon_m}-\frac{1}{\varepsilon_c}\right)\lambda_{j}
+\delta\lambda_{j,1}+o(\delta)+\mathcal{O}(\omega^2\delta)+\mathcal{O}(\omega^2)}.
\end{equation}
In (\ref{3.48}), $\lambda_{j}$ and $\varphi_{j}$ are respectively the eigenvalues and eigenfunctions of the operator $\mathcal{K}_{0}^*$. Moreover, we denote
\[
\langle\cdot,\cdot\rangle_{\mathcal{H}^*(S^c)}:=\langle \cdot,\cdot\rangle_{\mathcal{H}^*(S^a)}+\langle \cdot,\cdot\rangle_{\mathcal{H}^*(S^b)}.
\]
\end{theorem}

\begin{proof}
By Lemma \ref{freasy01} and using the eigenvalue and eigenfunction expansion formula in Lemma \ref{le:0102}, one can deduce (\ref{3.48}) by straightforward though tedious calculations.
\end{proof}

Next, we present the asymptotic expansion of the scattering field in \eqref{2.8} associated with the nanorod $D_\delta$.

\begin{theorem}\label{th:0101}
Let $u$ be the solution to (\ref{2.8}). Then, under conditions (C1) and (C2), the scattering field can be presented as
\begin{align}\label{3.53}
u^s(x)=&\sum_{j\in J}\frac{\mathrm{i}\omega\sqrt{\mu_m\varepsilon_m}a_{j,\delta}^{-1}\langle d\cdot\nu,\varphi_{j,\delta}\rangle\hat{\mathcal{S}}_{S^c}[\widetilde{\varphi}_j](x)\delta^2+o(\omega\delta^2)}
{\lambda\left(\frac{\varepsilon_m}{\varepsilon_c}\right)-\lambda_{j}
+\delta\left(\frac{1}{\varepsilon_c}-\frac{1}{\varepsilon_m}\right)^{-1}\lambda_{j,1}+o(\delta)+\mathcal{O}(\omega^2\delta)+\mathcal{O}(\omega^2)}\nonumber\\
&\ \ \ \ +\mathcal{O}(\omega\delta^3\|\tilde\psi\|_{\partial D})+\mathcal{O}(\omega\delta),
\end{align}
where
\beq
\lambda(t)=\frac{t+1}{2(t-1)},
\eeq
and the operator $\hat{\mathcal{S}}_{\Gamma}$, for any surface $\Gamma$, is defined by
\begin{equation}
\hat{\mathcal{S}}_{\Gamma}[\widetilde{\varphi}](x)=\int_{\Gamma}G^0(x-z_{\tilde{y}})\widetilde{\varphi}(\tilde{y})\, d\tilde{\sigma}(\tilde{y}),\label{eq:def01}
\end{equation}
\end{theorem}

\begin{proof}
First, by the Taylor expansion around $y=z_y$ for the Green function $G^{k_m}(x-y)$, and using (\ref{3.48}), we have that
\begin{align}\label{3.54}
u^s(x)&=\mathcal{S}_{D_\delta}^{k_m}[\psi](x)=\int_{\partial D_\delta}G^{k_m}(x-z_{y})\psi(y)d\sigma(y)\nonumber\\
&\ \ \ \ \ \ -\int_{\partial D_\delta}\nabla G^{k_m}(x-z_{y})\cdot(y-z_{y})\psi(y)d\sigma(y)+\mathcal{O}(\delta^3\|\tilde\psi\|_{\partial D}+\omega\delta).
\end{align}
Furthermore, noticing that
\begin{align*}
\partial D_\delta=S^a\cup\left(S^f\cap\overline{\iota_{1,\delta}(P_0)}\right)\cup\left(S^f\setminus\overline{\iota_{1,\delta}(P_0)\cup\iota_{1,\delta}(Q_0)}\right)
\cup\left(S^f\cap\overline{\iota_{1,\delta}(Q_0)}\right)\cup S^b,
\end{align*}
and then substituting (\ref{3.48}) into (\ref{3.54}) yield
\begin{align*}
u^s(x)&=\sum_{j\in J}\frac{\mathrm{i}\omega\sqrt{\mu_m\varepsilon_m}(1/\varepsilon_c-1/\varepsilon_m)a_{j,\delta}^{-1}\langle d\cdot\nu,\varphi_{j,\delta}\rangle\hat{\mathcal{S}}_{S^c}[\widetilde{\varphi}_j](x)+o(\omega\delta)+\mathcal{O}(\omega^2)}
{\frac{1}{2}\left(\frac{1}{\varepsilon_m}+\frac{1}{\varepsilon_c}\right)+\left(\frac{1}{\varepsilon_m}-\frac{1}{\varepsilon_c}\right)\lambda_{j}
+\delta\lambda_{j,1}+o(\delta)+\mathcal{O}(\omega^2\delta)+\mathcal{O}(\omega^2)}\delta^2\\
&\ \ \ +\sum_{j\in J}\frac{o(\omega\delta)}
{\frac{1}{2}\left(\frac{1}{\varepsilon_m}+\frac{1}{\varepsilon_c}\right)+\left(\frac{1}{\varepsilon_m}-\frac{1}{\varepsilon_c}\right)\lambda_{j}
+\delta\lambda_{j,1}+o(\delta)+\mathcal{O}(\omega^2\delta)+\mathcal{O}(\omega^2)}\delta\\
&\ \ \ -\sum_{j\in J}\frac{\mathrm{i}\omega\sqrt{\mu_m\varepsilon_m}(1/\varepsilon_c-1/\varepsilon_m)a_{j,\delta}^{-1}\langle d\cdot\nu,\varphi_{j,\delta}\rangle\int_{S^c}\frac{\partial G^0}{\partial\nu(\tilde{y})}(x-z_{\tilde{y}})\widetilde{\varphi}_j(\tilde{y})d\tilde{\sigma}(\tilde{y})}
{\frac{1}{2}\left(\frac{1}{\varepsilon_m}+\frac{1}{\varepsilon_c}\right)+\left(\frac{1}{\varepsilon_m}-\frac{1}{\varepsilon_c}\right)\lambda_{j}
+\delta\lambda_{j,1}+o(\delta)+\mathcal{O}(\omega^2\delta)+\mathcal{O}(\omega^2)}\delta^3\\
&\ \ \ -\sum_{j\in J}\frac{o(\omega\delta)}
{\frac{1}{2}\left(\frac{1}{\varepsilon_m}+\frac{1}{\varepsilon_c}\right)+\left(\frac{1}{\varepsilon_m}-\frac{1}{\varepsilon_c}\right)\lambda_{j}
+\delta\lambda_{j,1}+o(\delta)+\mathcal{O}(\omega^2\delta)+\mathcal{O}(\omega^2)}\delta^2\\
&\ \ \ +\mathcal{O}(\omega\delta)+\mathcal{O}(\delta^3\|\tilde\psi\|_{\partial D}),
\end{align*}
which proves the expansion formula (\ref{3.53}).

The proof is complete.
\end{proof}

\subsection{Quantitative analysis of the scattering wave field}\label{sect:2.4}

Based on the asymptotic expansion formula of the scattering field $u^s$ in Theorem \ref{th:0101}, we can perform some quantitative analysis of the scattering field associated with the nanorod $D_\delta$.

We first consider the case that the nanorod is straight, namely $\Gamma_0$ is a straight line with length $L$. Moreover, we assume that $S^a$ and $S^b$ are two semi-spheres.  Then using the definition in \eqnref{eq:def01}, the scattering field $u^s$ in \eqnref{3.53} is given by
\beq\label{eq:resca01}
u^s(x)=\mathrm{i}\omega\delta^2 p(x)\sum_{j\in J}\kappa_{j,\delta}+o(\delta^2\|\tilde\psi\|_{\partial D})+\mathcal{O}(\omega\delta),
\eeq
where $\kappa_{j,\delta}$ is defined by
\[
\kappa_{j,\delta}:=\frac{\sqrt{\mu_m\varepsilon_m}a_{j,\delta}^{-1}\langle d\cdot\nu,\varphi_{j,\delta}\rangle}
{\lambda\left(\frac{\varepsilon_m}{\varepsilon_c}\right)-\lambda_{j}+\delta\left(\frac{1}{\varepsilon_c}-\frac{1}{\varepsilon_m}\right)^{-1}\lambda_{j,1}
}\int_{S^a}\tilde\varphi_{j},
\]
and the function $p(x)$ has the form
\beq\label{eq:resca02}
p(x):=\frac{1}{|x-P_0|}+\frac{1}{|x-Q_0|}.
\eeq
We mention that in order to get \eqnref{eq:resca01}, we used the fact that
$$
\int_{S^a}\tilde\varphi_j=\int_{S^b}\tilde\varphi_j,
$$
which can be verified by using \eqnref{eq:K0def01}, \eqnref{eq:eigenasy00} and the fact that $S^a$ and $S^b$ are two semi-sphere with the same radius.
The amplitude of the scattering field $u^s$ is mainly determined by the function $p(x)$ in \eqnref{eq:resca02}.
By straightforward computations, one has
\beq
p(x)=\left\{
\begin{array}{ll}
\sqrt{(L/2-l(x))^2+\delta^2}+\sqrt{(L/2+l(x))^2+\delta^2}, & x\in S_\delta^f,\\
\delta+ \sqrt{\delta^2+L^2-2(x-P_0, P_0-Q_0)}   & x\in S_\delta^a,\\
\delta+ \sqrt{\delta^2+L^2-2(x-Q_0, Q_0-P_0)}   & x\in S_\delta^b,
\end{array}
\right.
\eeq
where
$$
l(x):=\Big|z_{x}-\frac{P_0+Q_0}{2}\Big|.
$$
Finally, by some elementary analysis, one can find that $p(x)$ attains its maximum $L+2\delta$ at the two ending points of $S_\delta^a$ and $S_\delta^b$, and attains its minimum $\sqrt{L^2+4\delta^2}$ on the centering parts (circular area) of $S_\delta^f$. The behaviuor is more obvious when $L$ is larger. Hence, we can conclude that the scattering field $u^s$ behaves stronger (in terms of the amplitude of the wave field) on the two end-parts of the nanorod $S_\delta^a$ and $S_\delta^b$ than that on the facade-part of the nanorod $S_\delta^f$.

\begin{figure}[htp]
\begin{center}
{\includegraphics[width=0.3\textwidth]{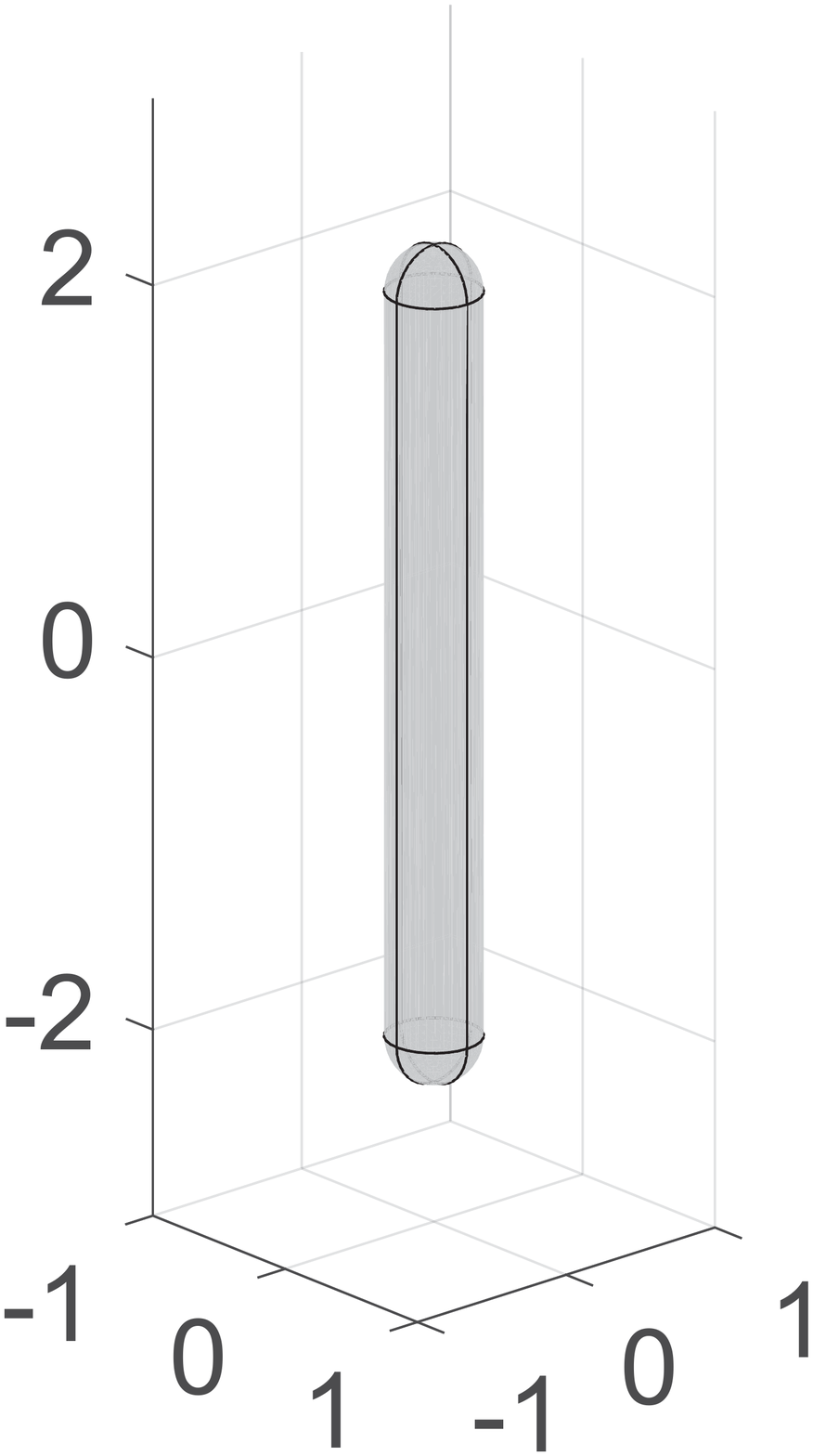}}
{\includegraphics[width=0.3\textwidth]{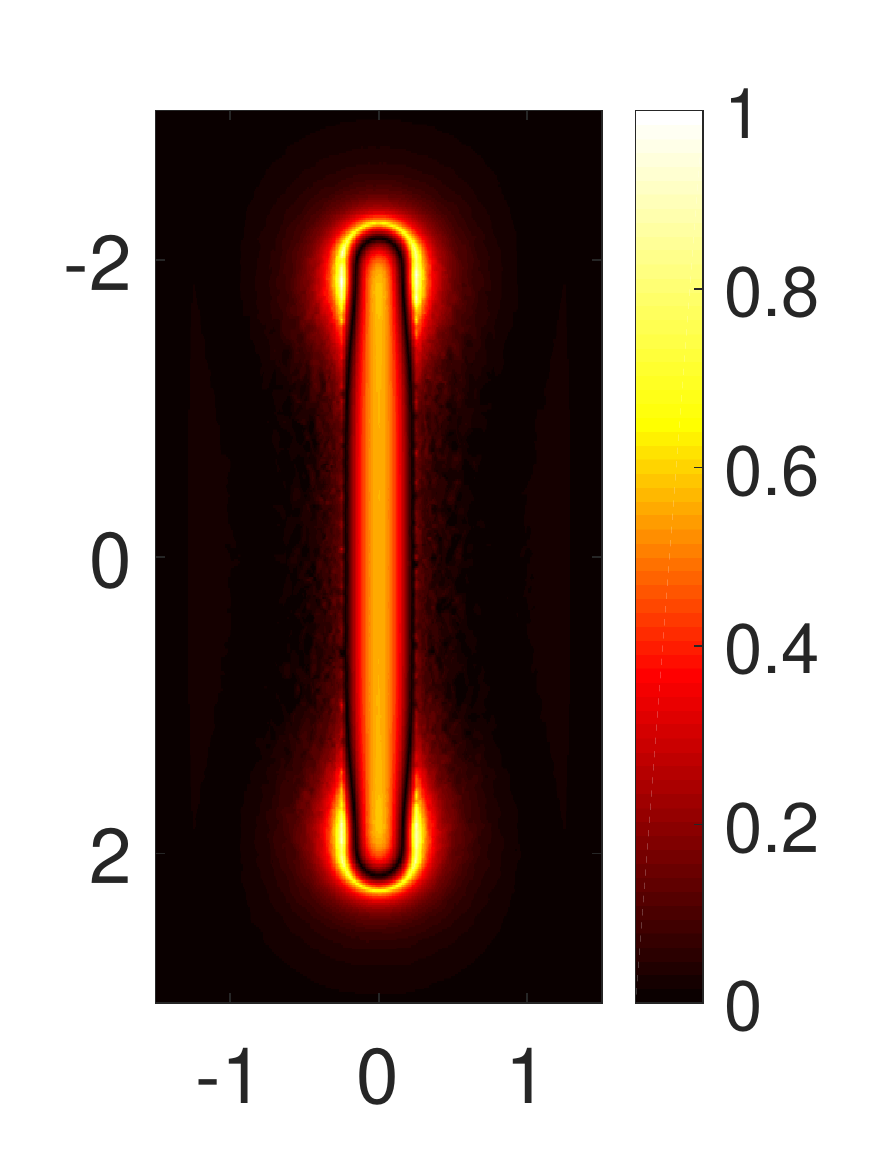}}
{\includegraphics[width=0.3\textwidth]{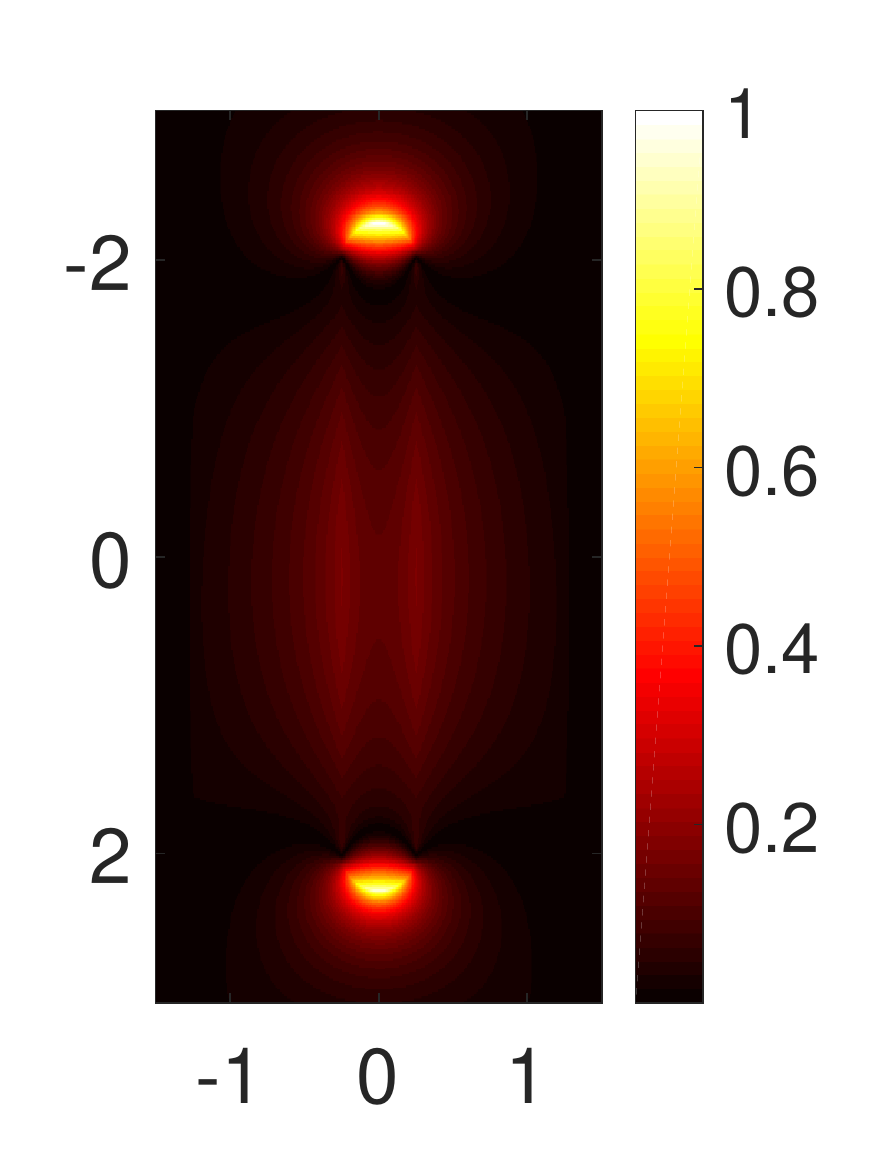}}
\end{center}
\caption{\label{fig:Plasmon01} Left: The geometry of a straight nanorod;  Middle: Slice plot of  the normalized scattering field $|\Re{u^s}|$ on the $(x_2,x_3)$-plane, where the incident direction is $d=(1,0,0)$;   Right: Slice plot of the normalized scattering field $|\Re{u^s}|$ on the $(x_2,x_3)$-plane, where the incident direction is $d=(0,0,1)$. The incident wave is: $u^i(x)=10^3 e^{\mathrm{i}k_m d\cdot x} $. }
\end{figure}

\begin{figure}[htp]
\begin{center}
{\includegraphics[width=0.3\textwidth]{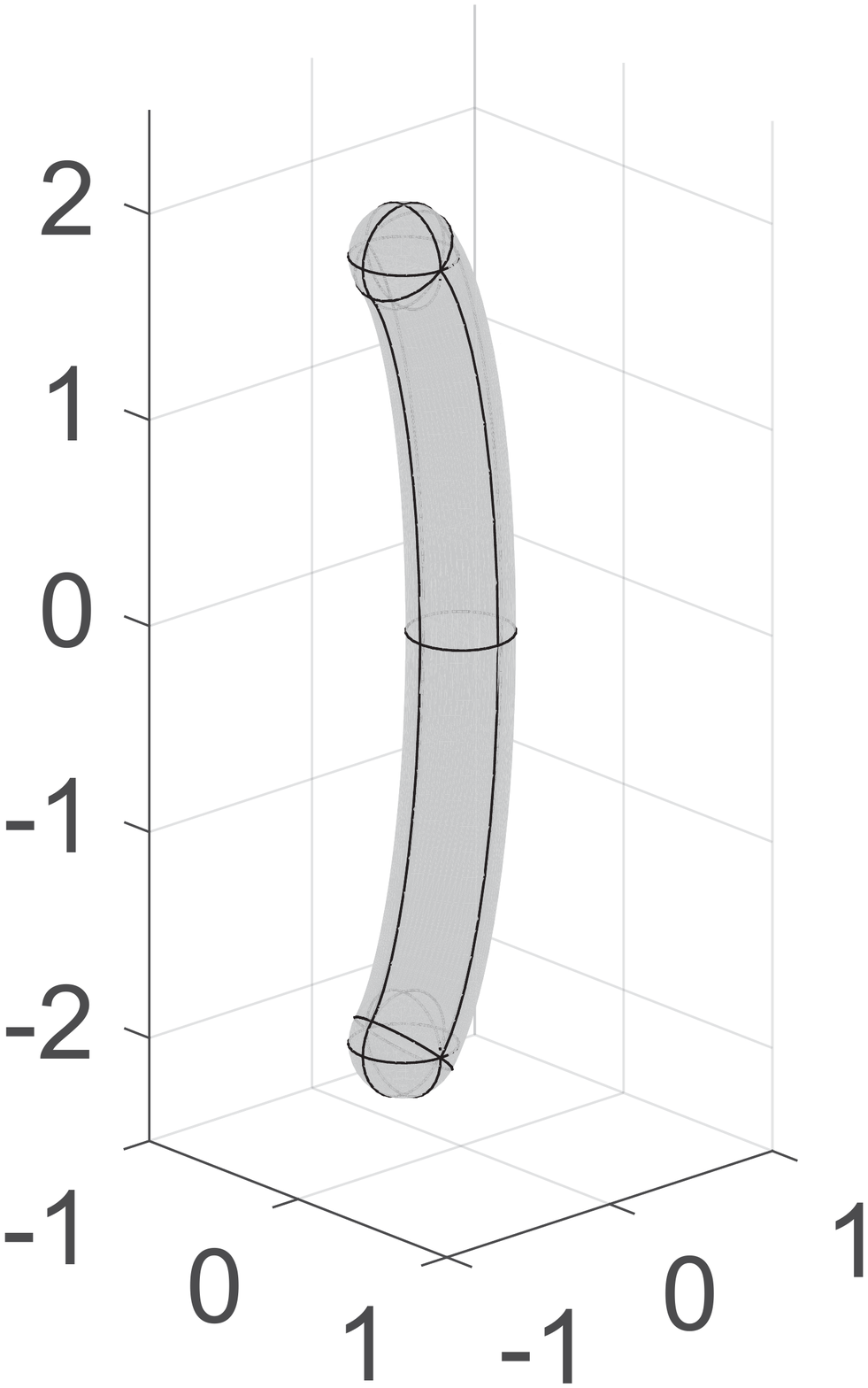}}
{\includegraphics[width=0.3\textwidth]{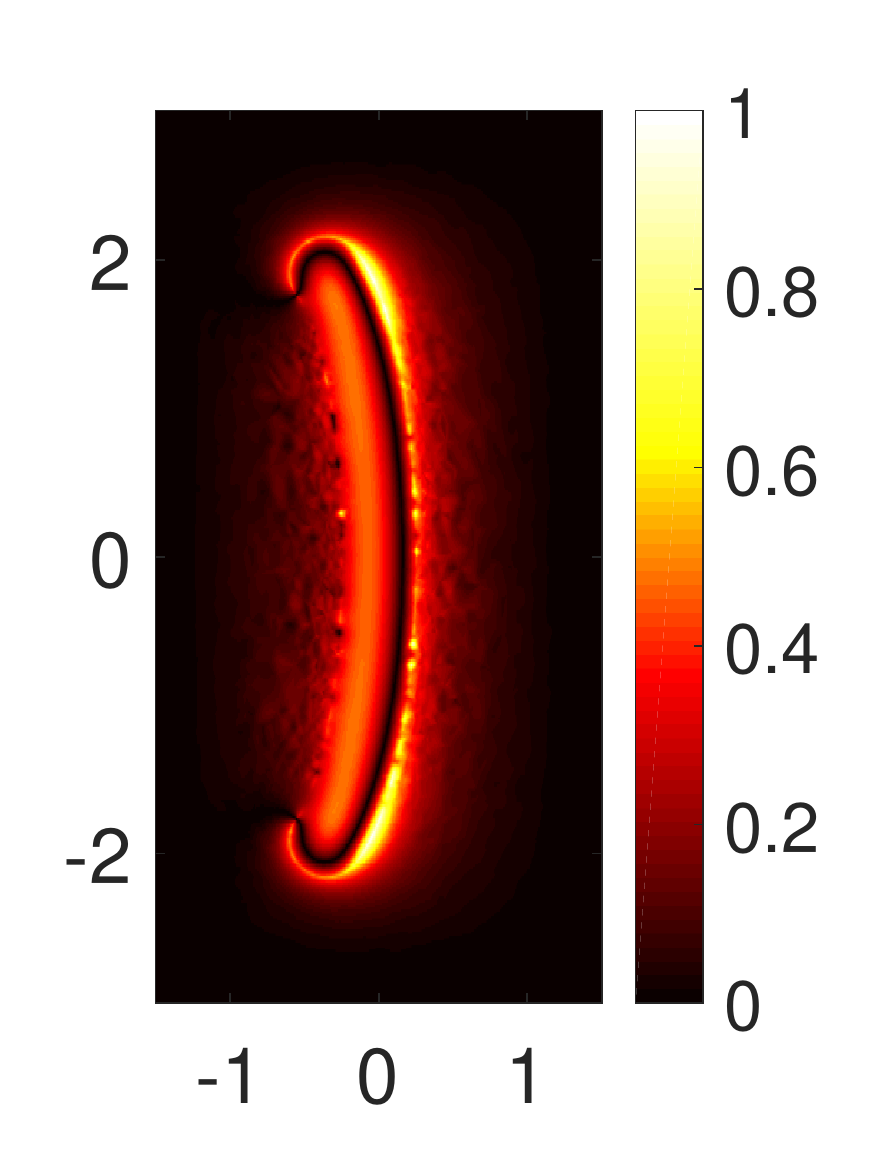}}
{\includegraphics[width=0.3\textwidth]{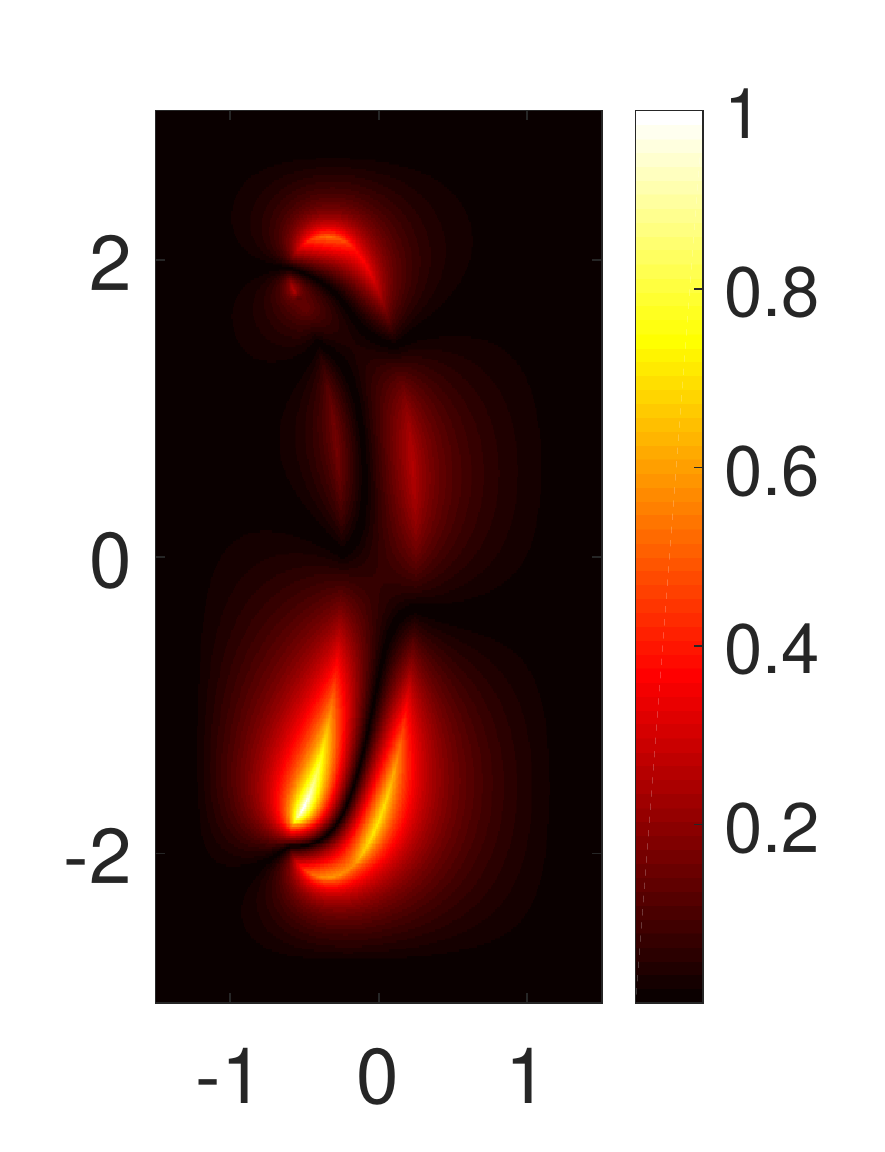}}
\end{center}
\caption{\label{fig:Plasmon02} Left: The geometry of a curved nanorod;  Middle: Slice plot of the normalized scattering field $|\Re{u^s}|$ on the $(x_2, x_3)$-plane, where the incident direction is $d=(1,0,0)$;   Right: Slice plot of the normalized scattering field $|\Re{u^s}|$ on the $(x_2, x_3)$-plane, where the incident direction is $d=(0,0,1)$. The incident wave is: $u^i(x)=10^3 e^{\mathrm{i}k_m d\cdot x} $.}
\end{figure}

For the general case with a generically curved nanorod, we believe the same quantitative behaviours should hold in a certain sense. In what follows, we present two numerical examples to illustrate such quantitative behaviours for the scattering field of the Helmholtz system \eqref{2.8} associated with a curved nanorod $D_\delta$; see Figures \ref{fig:Plasmon01} and \ref{fig:Plasmon02}. In the sequel, we adopt the coordinate notation $x=(x_j)_{j=1}^3\in\mathbb{R}^3$. In Fig.~\ref{fig:Plasmon01}, $D_\delta$ is generated by a straight $\Gamma_0$ (along the $x_3$-axis) of length 4 and two end-caps being two semi-balls of radius $\delta$; whereas in Fig.~\ref{fig:Plasmon02}, $D_\delta$ is generated by a curved $\Gamma_0$ and two end-caps being two semi-balls of radius $\delta$, where the parametrization of $\Gamma_0$ is $x(t)=(x_j(t))_{j=1}^3:$
\begin{equation}\label{eq:para1}
  \begin{cases}
  x_1(t)=0,\\
  x_2(t)=\frac{1}{2}(\cos(t)-1),\\
  x_3(t)=2\sin(t), \quad t\in[-\frac{\pi}{2}+\frac{3}{10}, \frac{\pi}{2}-\frac{3}{10} ].
  \end{cases}
\end{equation}
The other key parameters are given as follows,
\begin{equation}\label{eq:para2}
  \delta=\frac{1}{2^2},\quad \omega=\delta^{1/3}, \quad \varepsilon_c=-1+\mathrm{i}\omega^4,\quad \varepsilon_m=1,\quad \mu_{D_\delta}\equiv 1.
\end{equation}
Figures~\ref{fig:Plasmon01} and \ref{fig:Plasmon02} respectively plot the absolute value of the real parts of numerically computed wave fields, namely $|\Re{u^s}|$, associated with the two nanorods described above. It is remarked that $\Re u$ and $\Re u^s$ are the physical fields. Here, in order to present a better display, we normalize the wave fields in the sense that the maximum value of $|\Re{u^s}|$ in both figures is $1$. It is evident that the wave field attains its maximum amplitude at the two-end parts of the nanorod. It is emphasized that the choice of the material configuration in \eqref{eq:para2}, in particular $\Re\varepsilon_c=-1$ is mainly based on numerical simplicity and convenience. One can pick the other choices, say e.g. $\Re \varepsilon_c=-2$, and would have similar numerical behaviours. The numerical simulations are not the main focus of this paper, and we only present a typical example for illustration. On the other hand, we would like to point out that the material configuration \eqref{eq:para2} is nearly resonant according to Definition \ref{de2} with $\eta_0$ being set sufficiently small. We shall consider this example again in our subsequent resonance study and provide more relevant discussions in Section~\ref{sect:4.3} in what follows.



%

\section{Resonance analysis of the exterior wave field}

Using the asymptotic result in Theorem~\ref{th:0101}, we proceed to analyze the plasmon resonance of the scattering system \eqref{2.8}. We first derive the gradient estimate of the scattering field $u^s$ outside the nanorod $D_\delta$.
%

\begin{lemma}\label{lem4.1}
Let $\psi=\psi_c+\langle\psi,\varphi_{0,\delta}\rangle\varphi_{0,\delta}$, where $\psi_c\in\mathcal{H}_0^*(\partial D)$. Then for the scattering solution of (\ref{2.8}), i.e., $u^s=\mathcal {S}_{D_\delta}^{k_m}[\psi]$, we have the following estimate
\begin{align}\label{4.1}
\left|\|\nabla u^s\|_{L^2(\mathbb{R}^3\setminus\overline{D_\delta})}^2-\|\psi_c\|^2\right|\lesssim\omega\left|\langle\psi,\varphi_{0,\delta}\rangle\right|^2.
\end{align}
\end{lemma}

\begin{proof}
Let $B_R$ be a  sufficiently large ball such that $\overline{D_\delta}\subset B_R$. By using the divergence theorem in $B_R\setminus\overline{D_\delta}$, the jump relation \eqref{2.12} and the Sommerfeld radiation condition, it follows that
\begin{align*}
\int_{B_R\setminus\overline{D_\delta}}|\nabla u^s|^2dx= &\bar{k}_m^2\int_{B_R\setminus\overline{D_\delta}}|u^s|^2dx-\int_{\partial D_\delta}u^s\overline{\frac{\partial u^s}{\partial\nu}\Big|+}d\sigma+\int_{\partial B_R}u^s\overline{\frac{\partial u^s}{\partial R}}d\sigma\\
=&\bar{k}_m^2\int_{B_R\setminus\overline{D_\delta}}|u^s|^2dx-\int_{\partial D_\delta}\mathcal {S}_{D_\delta}^{k_m}[\psi]\overline{\left(\frac{1}{2}\mathcal{I}+(\mathcal{K}_{D_\delta}^{k_m})^*\right)[\psi]}d\sigma\\
&+\int_{\partial B_R}u^s\cdot\overline{ik_mu^s+\mathcal{O}(R^{-2})}d\sigma.
\end{align*}
From the expansion formulas (\ref{8.1}) and (\ref{8.6}), we obtain
\begin{align}
\left|\int_{\partial D_\delta}\mathcal {S}_{D_\delta}^{k_m}[\psi]\overline{\left(\frac{1}{2}\mathcal{I}+(\mathcal{K}_{D_\delta}^{k_m})^*\right)[\psi]}d\sigma\right|
\leq\left|\int_{\partial D_\delta}\mathcal {S}_{D_\delta}[\psi]\overline{\left(\frac{1}{2}\mathcal{I}+\mathcal{K}_{D_\delta}^*\right)[\psi]}d\sigma\right|+\left|E\right|,
\end{align}
where
\begin{align}
E=\int_{\partial D_\delta}\mathcal{S}_{D_\delta}[\psi]\overline{\sum_{j=1}^\infty k_m^{j}\mathcal{K}_{D_\delta,j}[\psi]}d\sigma
+\int_{\partial D_\delta}\sum_{j=1}^\infty k_m^{j}\mathcal
{S}_{D_\delta,j}[\psi]\overline{\left(\frac{1}{2}\mathcal
{I}+(\mathcal{K}_{D_\delta}^{k_m})^*\right)[\psi]}d\sigma\nonumber.
\end{align}
Owing to $\omega$ is sufficiently small, by Cauchy's inequality, it is easy to see that $|E|\lesssim k_m\|\psi\|^2\lesssim \omega\|\psi\|^2$. Next, we estimate $\left|\int_{\partial D_\delta}\mathcal {S}_{D_\delta}[\psi]\overline{\left(\frac{1}{2}\mathcal{I}+\mathcal{K}_{D_\delta}^*\right)[\psi]}d\sigma\right|$.
Since $\mathcal{K}_{D_\delta}^*[\varphi_{0,\delta}]=\frac{1}{2}\varphi_{0,\delta}$ and $\mathcal{S}_{D_\delta}[\varphi_{0,\delta}]=0$, for $\psi=\psi_c+\langle\psi,\varphi_{0,\delta}\rangle\varphi_{0,\delta}$, it implies that
\begin{align*}
&\int_{\partial D_\delta}\mathcal {S}_{D_\delta}[\psi]\overline{\left(\frac{1}{2}\mathcal{I}+\mathcal{K}_{D_\delta}^*\right)[\psi]}d\sigma\\
=&\int_{\partial D_\delta}\mathcal {S}_{D_\delta}[\psi_c]\overline{\left(\langle\psi,\varphi_{0,\delta}\rangle\varphi_{0,\delta}+\left(\frac{1}{2}\mathcal{I}+\mathcal{K}_{D_\delta}^*\right)[\psi_c]\right)}d\sigma\\
=&\int_{\partial D_\delta}\mathcal{S}_{D_\delta}\left[\sum_{j=1}^\infty\langle\psi_c,\varphi_{j,\delta}\rangle\varphi_{j,\delta}\right]\overline{\left(\langle\psi,\varphi_{0,\delta}\rangle\varphi_{0,\delta}
+\left(\frac{1}{2}\mathcal{I}+\mathcal{K}_{D_\delta}^*\right)\left[\sum_{l=1}^\infty\langle\psi_c,\varphi_{l,\delta}\rangle\varphi_{l,\delta}\right]\right)}d\sigma\\
=&\int_{\partial D_\delta}\sum_{j=1}^\infty\langle\psi_c,\varphi_{j,\delta}\rangle\mathcal{S}_{D_\delta}[\varphi_{j,\delta}]\overline{\left(\langle\psi,\varphi_{0,\delta}\rangle\varphi_{0,\delta}
+\sum_{l=1}^\infty\langle\psi_c,\varphi_{l,\delta}\rangle\left(\frac{1}{2}+\lambda_{l,\delta}\right)\varphi_{l,\delta}\right)}d\sigma\\
=&\sum_{j=1}^\infty\overline{\langle\psi,\varphi_{0,\delta}\rangle}\langle\psi_c,\varphi_{j,\delta}\rangle\int_{\partial D_\delta}\mathcal{S}_{D_\delta}[\varphi_{j,\delta}]\overline{\varphi_{0,\delta}}d\sigma\\
&\quad +\sum_{j,l=1}^\infty\left(\frac{1}{2}+\lambda_{l,\delta}\right)\overline{\langle\psi_c,\varphi_{l,\delta}\rangle}\langle\psi_c,\varphi_{j,\delta}\rangle\int_{\partial D_\delta}\mathcal{S}_{D_\delta}[\varphi_{j,\delta}]\overline{\varphi_{l,\delta}}d\sigma.
\end{align*}
Noting that
\begin{align*}
\int_{\partial D_\delta}\mathcal{S}_{D_\delta}[\varphi_{j,\delta}]\overline{\varphi_{l,\delta}}d\sigma=-\langle\varphi_{l,\delta},\varphi_{j,\delta}\rangle=
\begin{cases}
1,\ \ \ l=j,\\
0,\ \ \ l\neq j,
\end{cases}
\end{align*}
it deduces that
\begin{align*}
\int_{\partial D_\delta}\mathcal {S}_{D_\delta}[\psi]\overline{\left(\frac{1}{2}\mathcal{I}+\mathcal{K}_{D_\delta}^*\right)[\psi]}d\sigma
=\sum_{j=1}^\infty\left(\frac{1}{2}+\lambda_{j,\delta}\right)\left|\langle\psi_c,\varphi_{j,\delta}\rangle\right|^2.
\end{align*}
Since $\lambda_{j,\delta}\in\left(-\frac{1}{2},\frac{1}{2}\right)$ $(j\geq1)$, we can obtain
\begin{align}\label{4.3}
\left|\int_{\partial D_\delta}\mathcal {S}_{D_\delta}[\psi]\overline{\left(\frac{1}{2}\mathcal{I}+\mathcal{K}_{D_\delta}^*\right)[\psi]}d\sigma
\right|\approx\|\psi_c\|^2.
\end{align}
Moreover, by Cauchy's inequality, it follows
\begin{align}\label{4.4}
\left|\int_{\partial B_R}u^s\cdot\overline{\mathrm{i}k_mu^s+\mathcal{O}(R^{-2})}d\sigma\right|\lesssim&\omega\|u^s\|_{L^2(\partial B_R)}^2+\int_{\partial B_R}|u^s\cdot\overline{\mathcal{O}(R^{-2})}|d\sigma\nonumber\\
\lesssim&\omega\|\psi\|^2+\mathcal{O}(R^{-1})\cdot\|\psi\|.
\end{align}
By combing (\ref{4.3}) and (\ref{4.4}), for $\omega\ll1$, we obtain
\begin{align*}\label{}
\|\nabla u^s\|_{L^2(B_R\setminus\overline{D_\delta})}^2\lesssim\|\psi_c\|^2+\omega\left|\langle\psi,\varphi_{0,\delta}\rangle\right|^2+\mathcal{O}(R^{-1})\cdot\|\psi\|.
\end{align*}
Similarly, by (\ref{4.3}) and (\ref{4.4}), we also deduce the inverse inequality as
\begin{align*}\label{}
\|\nabla u^s\|_{L^2(B_R\setminus\overline{D_\delta})}^2\gtrsim\|\psi_c\|^2-\omega\left|\langle\psi,\varphi_{0,\delta}\rangle\right|^2-\mathcal{O}(R^{-1})\cdot\|\psi\|.
\end{align*}
Hence, by letting $R\rightarrow\infty$, we see that the estimate (\ref{4.1}) holds.

The proof is complete.
\end{proof}

Before proceeding with the gradient analysis of the scattering field $u^s$ outside the nanorod $D_\delta$, we consider the parameter choice of the permittivity with an imaginary part. In fact, in real applications, nano-metal materials always contain losses, which are reflected in the imaginary part of a complex electric permittivity $\varepsilon_c$. As shall be shown, like the frequency $\omega$ and the size $\delta$, the lossy parameter also plays a key role in the plasmon resonance of the scattering field of the nanorod. Let $\theta=\Re\left(\frac{1}{\varepsilon_c}\right)$, $\rho=\Im\left(\frac{1}{\varepsilon_c}\right)<0$. Then
$\tau_{j,\delta}$ given by (\ref{3.4}) can be written as
\begin{align}\label{4.5}
\tau_{j,\delta}=\frac{1}{2}(\theta+\varepsilon_m^{-1})-(\theta-\varepsilon_m^{-1})\lambda_{j,\delta}+\rho(\frac{1}{2}-\lambda_{j,\delta})\mathrm{i}.
\end{align}
Next, by considering the principal equation $\mathcal{A}_{D_\delta,0}[\psi_0]=f_0$, where $\mathcal{A}_{D_\delta,0}$ is defined by (\ref{2.20}) and $\psi_0, f_0\in \mathcal{H}^*(\partial D_\delta)$, similar to (\ref{3.7}), applying the eigenfunction expansion, it follows that
\begin{align}\label{4.6}
\psi_0=\mathcal {A}_{D_\delta,0}^{-1}[f_0]=\sum_{j=0}^\infty\frac{a_{j\delta}^{-1}\langle f_0,\varphi_{j,\delta}\rangle}{\tau_{j,\delta}}\varphi_{j,\delta}.
\end{align}
\begin{lemma}\label{lem4.2}
Under conditions (C1) and (C2), $\psi_0$ is given by (\ref{4.6}) and has the decomposition $\psi_0=\psi_{0,c}+c\varphi_0$, ($\psi_{0,c}\in\mathcal{H}_0^*(\partial D_\delta)$, $c$ is a constant). Then, for sufficiently small $|\rho|$, it holds that
\begin{enumerate}

\item[(1)] $\|\mathcal{A}_{D_\delta,0}^{-1}\|_{\mathcal{L}(\mathcal{H}^*(\partial D_\delta),\mathcal{H}^*(\partial D_\delta))}\lesssim|\rho|^{-1}$.

\item[(2)] If $\frac{1}{2}\frac{\theta+\varepsilon_m^{-1}}{\theta-\varepsilon_m^{-1}}\neq\lambda_{j,\delta}, (j\geq0)$, then $\|\mathcal{A}_{D_\delta,0}^{-1}\|_{\mathcal{L}(\mathcal{H}^*(\partial D_\delta),\mathcal{H}^*(\partial D_\delta))}\lesssim C$ for some positive constant $C$.

\item[(3)] If $\frac{1}{2}\frac{\theta+\varepsilon_m^{-1}}{\theta-\varepsilon_m^{-1}}=\lambda_{j,\delta}$ for some $j\geq1$, then $\|\psi_{0,c}\|\gtrsim|\rho|^{-1}\left|\langle f_0,\varphi_{j,\delta}\rangle\right|$.
\end{enumerate}
\end{lemma}

\begin{proof}
(1) For $j\neq0$, since $\left|\tau_{j,\delta}^{-1}\right|\lesssim\frac{1}{|\rho|(\frac{1}{2}-\lambda_{j,\delta})}\lesssim|\rho|^{-1}$, it follows that
\begin{equation}\label{}
\|\psi_0\|^2\lesssim|\rho|^{-2}\sum_{j=0}^\infty\left|\langle f_0,\varphi_{j,\delta}\rangle\right|^2\lesssim|\rho|^{-2}\|f_0\|^2.
\end{equation}
Hence, $\|\mathcal{A}_{D_\delta,0}^{-1}\|_{\mathcal{L}(\mathcal{H}^*(\partial D_\delta),\mathcal{H}^*(\partial D_\delta))}\lesssim|\rho|^{-1}$.

\medskip

(2) If $\frac{1}{2}\frac{\theta+\varepsilon_m^{-1}}{\theta-\varepsilon_m^{-1}}\neq\lambda_{j,\delta}$, then, for $j\geq0$, we have $\left|\frac{1}{2}\frac{\theta+\varepsilon_m^{-1}}{\theta-\varepsilon_m^{-1}}-\lambda_{j,\delta}\right|\geq c_0$, where $c_0$
is a positive constant. Therefore, $\left|\tau_{j,\delta}^{-1}\right|\lesssim1$ and $\|\psi_0\|^2\lesssim\|f_0\|^2$, i.e.,
$\|\mathcal{A}_{D_\delta,0}^{-1}\|_{\mathcal{L}(\mathcal{H}^*(\partial D_\delta),\mathcal{H}^*(\partial D_\delta))}\lesssim C$.

\medskip

(3) When $\frac{1}{2}\frac{\theta+\varepsilon_m^{-1}}{\theta-\varepsilon_m^{-1}}=\lambda_{j,\delta}$ for some $j\geq1$, by (\ref{4.5}), one has
$\tau_{j,\delta}=\rho(\frac{1}{2}-\lambda_{j,\delta})\mathbb{}\mathrm{i}$, it then follows that
\begin{equation}\label{}
\|\psi_{0,c}\|\gtrsim\left|\langle\psi_{0,c},\varphi_{j,\delta}\rangle\right|\gtrsim\frac{\left|\langle f_0,\varphi_{j,\delta}\rangle\right|}{\left|\rho(\frac{1}{2}-\lambda_{j,\delta})\right|}\gtrsim|\rho|^{-1}\left|\langle f_0,\varphi_{j,\delta}\rangle\right|,
\end{equation}
which completes the proof.

\end{proof}

With Lemma \ref{lem4.2}, we can establish the following key result for estimating the gradient of the scattering field $u^s$, which provide resonant and non-resonant conditions for the scattering system \eqref{2.8} associated with the nanorod $D_\delta$ according to criterion \eqref{prdf01} in Definition~\ref{depr}.

\begin{theorem}\label{thm4.3}
Let $u^s$ be the scattering solution of (\ref{2.8}). Suppose that $|\rho|^{-1}\omega^2\delta\leq c_1$ for a sufficiently small $c_1$, then under conditions (C1) and (C2), we have the following results.

\begin{enumerate}

\item[(1)] If $\frac{1}{2}\frac{\theta+\varepsilon_m^{-1}}{\theta-\varepsilon_m^{-1}}\neq\lambda_{j,\delta}$ for any $j\geq0$, there exists a constant $C$ independent of $\delta$ such that
\begin{equation}
\left\|\nabla u^s\right\|_{L^2(\mathbb{R}^3\setminus\overline{D_\delta})}\leq C.
\end{equation}

\item[(2)] If $\frac{1}{2}\frac{\theta+\varepsilon_m^{-1}}{\theta-\varepsilon_m^{-1}}=\lambda_{j,\delta}$ and $\langle d\cdot\nu,\varphi_{j,\delta}\rangle\neq0$ for some $j\geq1$, it holds
\begin{align}\label{bl1}
\left\|\nabla u^s\right\|_{L^2(\mathbb{R}^3\setminus\overline{D_\delta})}\gtrsim\mathcal{O}(|\rho|^{-1}\omega)+o(|\rho|^{-1}\omega\delta)
+\mathcal{O}(|\rho|^{-1}\omega^2)+\mathcal{O}(\omega^{\frac{1}{2}}).
\end{align}
Furthermore, assuming that  $|\rho|=o(\omega)$ and $|\rho|^{-1}\omega^2\delta\rightarrow 0$ as $\omega\rightarrow 0$, then it follows that
\begin{align}\label{bl2}
\left\|\nabla u^s\right\|_{L^2(\mathbb{R}^3\setminus\overline{D_\delta})}\rightarrow\infty\quad \mbox{as}\ \ \omega\rightarrow 0.
\end{align}
\end{enumerate}
\end{theorem}

\begin{proof}
(1) From (\ref{3.8}), it deduces that
\begin{equation}\label{}
\mathcal {A}_{D_\delta}(\omega)=\mathcal{A}_{D_\delta,0}\left(\mathcal{I}+\mathcal{A}_{D_\delta,0}^{-1}\left(\omega^{2}\delta\widehat{\mathcal
{A}}_{D_\delta,2}+o(\omega^2\delta)+\mathcal{O}(\omega^3)\right)\right),
\end{equation}
and then
\begin{equation}\label{}
\psi=\left(\mathcal{I}+\mathcal{A}_{D_\delta,0}^{-1}\left(\omega^{2}\delta\widehat{\mathcal
{A}}_{D_\delta,2}+o(\omega^2\delta)+\mathcal{O}(\omega^3)\right)\right)^{-1}\mathcal{A}_{D_\delta,0}^{-1}[f].
\end{equation}
By Lemma \ref{lem4.2} (1), one finds
\begin{equation}\label{}
\left\|\mathcal{A}_{D_\delta,0}^{-1}\left(\omega^{2}\delta\widehat{\mathcal
{A}}_{D_\delta,2}+o(\omega^2\delta)+\mathcal{O}(\omega^3)\right)\right\|_{\mathcal{L}(\mathcal{H}^*(\partial D_\delta),\mathcal{H}^*(\partial D_\delta))}\lesssim|\rho|^{-1}\omega^2\delta.
\end{equation}
Hence, it follows that
\begin{align}
\left\|\psi-\psi_0\right\|=&\left\|\left(\mathcal{I}+\mathcal{A}_{D_\delta,0}^{-1}\left(\omega^{2}\delta\widehat{\mathcal
{A}}_{D_\delta,2}+o(\omega^2\delta)+\mathcal{O}(\omega^3)\right)\right)^{-1}\mathcal{A}_{D_\delta,0}^{-1}[f]-\mathcal{A}_{D_\delta,0}^{-1}[f]\right\|\nonumber\\
\lesssim&|\rho|^{-1}\omega^2\delta\left\|\mathcal{A}_{D_\delta,0}^{-1}[f]\right\|\nonumber\\
=&|\rho|^{-1}\omega^2\delta\left\|\psi_0\right\|.\label{4.15}
\end{align}
If $\frac{1}{2}\frac{\theta+\varepsilon_m^{-1}}{\theta-\varepsilon_m^{-1}}\neq\lambda_{j,\delta}$, by Lemma \ref{lem4.2} (2), we obtain
\begin{align*}
\left\|\psi\right\|\lesssim(1+|\rho|^{-1}\omega^2\delta)\left\|\psi_0\right\|
=(1+|\rho|^{-1}\omega^2\delta)\left\|\mathcal{A}_{D_\delta,0}^{-1}[f]\right\|\lesssim(1+c_1)\left\|f\right\|.
\end{align*}
Then, from Lemma \ref{lem4.1}, it yields
\begin{align*}
\left\|\nabla u^s\right\|_{L^2(\mathbb{R}^3\setminus\overline{D_\delta})}^2\lesssim\left\|\psi\right\|^2\lesssim C.
\end{align*}

(2) If $\frac{1}{2}\frac{\theta+\varepsilon_m^{-1}}{\theta-\varepsilon_m^{-1}}=\lambda_{j,\delta}$, ($j\geq1$), then, by using Lemma \ref{lem4.2} (3), we have that
\begin{align}\label{4.16}
\|\psi_{0,c}\|\gtrsim|\rho|^{-1}\left|\langle f,\varphi_{j,\delta}\rangle\right|.
\end{align}
Moreover, from (\ref{4.15}), it implies
\begin{align}\label{4.17}
|\rho|^{-1}\omega^2\delta\left\|\psi_0\right\|\gtrsim\left\|\psi-\psi_0\right\|\gtrsim\left\|\psi_c-\psi_{0,c}\right\|
\gtrsim\left\|\psi_{0,c}\right\|-\left\|\psi_c\right\|.
\end{align}

Combining now (\ref{4.16}) and (\ref{4.17}), and noticing that $|\rho|^{-1}\omega^2\delta\leq c_1$ for a sufficiently small $c_1$, we have
\begin{align}\label{}
\left\|\psi_c\right\|\gtrsim\left\|\psi_{0,c}\right\|-|\rho|^{-1}\omega^2\delta\left\|\psi_0\right\|\gtrsim|\rho|^{-1}\left|\langle f,\varphi_{j,\delta}\rangle\right|.
\end{align}
Therefore, we obtain from Lemma \ref{lem4.1} that
\begin{align*}
\left\|\nabla u^s\right\|_{L^2(\mathbb{R}^3\setminus\overline{D_\delta})}^2&\gtrsim
\left\|\psi_c\right\|^2-\omega\left|\langle \psi,\varphi_{0,\delta}\rangle\right|^2\\
&\gtrsim|\rho|^{-2}\left|\langle f,\varphi_{j,\delta}\rangle\right|^2-\omega\left|\langle \psi,\varphi_{0,\delta}\rangle\right|^2.
\end{align*}
Furthermore, applying Theorem 3.9, it follows that
\begin{align*}
\langle f,\varphi_{j,\delta}\rangle=\mathrm{i}\omega\sqrt{\mu_m\varepsilon_m}(1/\varepsilon_c-1/\varepsilon_m)a_{j,\delta}^{-1}\langle d\cdot\nu,\varphi_{j,\delta}\rangle\widetilde\varphi_{j}(\tilde{x})+o(\omega\delta)+\mathcal{O}(\omega^2).
\end{align*}
Thus, by noting that $|\langle \psi,\varphi_{0,\delta}\rangle|$ is bounded, (\ref{bl1}) holds and then by choosing $|\rho|=o(\omega)$ one has \eqnref{bl2}.

The proof is complete.
\end{proof}

\begin{remark}\label{re:0101}
It is noted that $\omega^{2}\delta=o(|\rho|)$ implies that $|\rho|^{-1}\omega^2\delta\rightarrow0$ as $\omega\rightarrow0$. From (2) of Theorem \ref{thm4.3}, one readily sees that when $\frac{1}{2}\frac{\theta+\varepsilon_m^{-1}}{\theta-\varepsilon_m^{-1}}=\lambda_{j,\delta}$, $\omega^{2}\delta=o(|\rho|)$ and $|\rho|=o(\omega)$, the gradient of the scattering field blows up $\omega\rightarrow 0$. According to \eqref{prdf01} in Definition \ref{depr}, we see that the plasmon resonance occurs. Notice that the last two conditions on the resonant material configuration are very flexible, and for example, one can take $|\rho|=\omega^s$, $s> 1$, $\delta=o(\omega^{s-2})$.
\end{remark}

\begin{remark}
Note that case (2) in Theorem \ref{4.3} is the resonance condition, i.e. $\frac{1}{2}\frac{\theta+\varepsilon_m^{-1}}{\theta-\varepsilon_m^{-1}}=\lambda_{j,\delta}$, $(j\geq1)$. If the lossy parameter of the nanorod $\Im(\varepsilon_c)\rightarrow0$, the resonance condition is consistent with $\tau_{j,\delta}=\frac{1}{2}\left(\frac{1}{\varepsilon_m}+\frac{1}{\varepsilon_c}\right)+\left(\frac{1}{\varepsilon_m}-\frac{1}{\varepsilon_c}\right)\lambda_{j,\delta}\rightarrow 0$
, which appeared in \eqnref{3.4} and \eqnref{3.7}.
\end{remark}

\section{Resonance analysis of the interior wave field}

\subsection{Asymptotics of the interior field in the quasi-static regime}
Owing to the relationship between the internal density $\phi$ and the scattering density $\psi$, we can deduce the asymptotic formula of the interior field with respect to $\omega$ as follows.

\begin{lemma}\label{interdens01}
In the quasi-static regime, under conditions (C1) and (C2), the interior field $u|_{D_\delta}$ has
the following representation
\begin{equation}\label{}
u(x)=\mathcal {S}_{D_\delta}^{k_c}[\phi](x),\ \ \ x\in D_\delta,
\end{equation}
where
\begin{equation}
\phi=\sum_{j\in J}\frac{\mathrm{i}\omega\sqrt{\mu_m\varepsilon_m}(1/\varepsilon_c-1/\varepsilon_m)a_{j,\delta}^{-1}\langle d\cdot\nu,\varphi_{j,\delta}\rangle\varphi_{j,\delta}+\mathcal{O}(\omega^2)}{\tau_{j,\delta}+\mathcal{O}(\omega^2)}+c\varphi_{0,\delta}+\mathcal{O}(\omega),
\end{equation}
and $c$ is a constant.
\end{lemma}

\begin{proof}
Substituting (\ref{8.1}), (\ref{8.3}) into (\ref{2.16}), and using the Taylor expansion $u^i=1+\mathrm{i}k_m d\cdot x +\mathcal{O}(\omega^2)$, we conclude
\begin{align}\label{relation}
\phi&=(\mathcal{S}_{D_\delta}^{-1}+k_c\mathcal{B}_{D_\delta,1}+k_c^2\mathcal{B}_{D_\delta,2}+\cdots)\left[(\mathcal{S}_{D_\delta}+\sum_{j=1}^\infty k_m^{j}\mathcal{S}_{D_\delta,j})[\psi]+1+\mathrm{i}k_m d\cdot x +\mathcal{O}(\omega^2)\right]\nonumber\\
&=\psi+\mathcal{S}_{D_\delta}^{-1}[1]+\mathcal{O}(\omega).
\end{align}
From (\ref{3.7}), and noticing $\mathcal{S}_{D_\delta}^{-1}[1]=c\varphi_{0,\delta}$ ($c$ is a constant), it follows that
\[
\phi=\sum_{j\in J}\frac{\mathrm{i}\omega\sqrt{\mu_m\varepsilon_m}(1/\varepsilon_c-1/\varepsilon_m)a_{j,\delta}^{-1}\langle d\cdot\nu,\varphi_{j,\delta}\rangle\varphi_{j,\delta}+\mathcal{O}(\omega^2)}{\tau_{j,\delta}+\mathcal{O}(\omega^2)}+c\varphi_{0,\delta}+\mathcal{O}(\omega),
\]
which readily completes the proof.
\end{proof}

In what follows, we shall make use of Lemma \ref{interdens01} to derive the asymptotic form of the internal energy.

\subsection{Asymptotics of the interior field with respect to the size parameter $\delta$}
\begin{theorem}\label{}
Let $u$ be the solution to (\ref{2.8}). Then, under conditions (C1) and (C2), for every $x\in D_\delta$, the interior field can be presented as
\begin{align}\label{interfield01}
u(x)=&\sum_{j\in J}\frac{\omega \mathrm{i}\sqrt{\frac{\mu_m}{\varepsilon_m}}a_{j,\delta}^{-1}\langle d\cdot\nu,\varphi_{j,\delta}\rangle\hat{\mathcal{S}}_{S^c}[\widetilde{\varphi}_j](x)\delta^2+o(\omega\delta^2)}
{\lambda\left(\frac{\varepsilon_m}{\varepsilon_c}\right)-\lambda_{j}
+\delta\left(\frac{1}{\varepsilon_c}-\frac{1}{\varepsilon_m}\right)^{-1}\lambda_{j,1}+o(\delta)+\mathcal{O}(\omega)}\nonumber\medskip\\
&\ \ \ \ +o(\delta)+\mathcal{O}(\delta^3\|\tilde\phi\|_{\partial D})+\mathcal{O}(\omega\delta),
\end{align}
where $a_{j,\delta}$, $\hat{\mathcal{S}}_{S^c}$, $\lambda_{j,1}$ and function $\lambda(t)$ are defined in Theorem \ref{th:0101}.
\end{theorem}
\begin{proof}
From Lemmas \ref{interdens01} and \ref{le:0102}, it follows that
\begin{align*}
\phi=\sum_{j\in J}\frac{\mathrm{i}\omega\sqrt{\mu_m\varepsilon_m}(1/\varepsilon_c-1/\varepsilon_m)a_{j,\delta}^{-1}\langle d\cdot\nu,\varphi_{j,\delta}\rangle\widetilde{\varphi}_{j}(\tilde{x})+o(\omega\delta)}
{\frac{1}{2}\left(\frac{1}{\varepsilon_m}+\frac{1}{\varepsilon_c}\right)+\left(\frac{1}{\varepsilon_m}-\frac{1}{\varepsilon_c}\right)\lambda_{j}
+\delta\lambda_{j,1}+o(\delta)+\mathcal{O}(\omega)}+c\widetilde{\varphi}_{0}(\tilde{x})+o(\delta)+\mathcal{O}(\omega).
\end{align*}
Substituting the above identity into the interior field formula $u=\mathcal{S}_{D_\delta}^{k_c}[\phi]$, and the by following a similar argument to that in the proof of Theorem \ref{th:0101}, one can prove
(\ref{interfield01}).
\end{proof}

We are now in a position to derive the energy estimate of the scattering system \eqref{2.8} within the curved nanorod $D_\delta$. Specifically, we shall show that under the same resonant condition in Theorem~\ref{thm4.3} (which characterises the resonance of the wave field outside the nanorod $D_\delta$), the internal energy of the nanorod blows up as well.

\begin{lemma}\label{interle01}
Let $\phi=\phi_c+\langle\phi,\varphi_{0,\delta}\rangle\varphi_{0,\delta}$, where $\phi_c\in\mathcal{H}_0^*(\partial D)$. Then for the interior solution of (\ref{2.8}) in $D_\delta$, i.e., $u=\mathcal {S}_{D_\delta}^{k_c}[\phi]$ in $D_\delta$, we have the following estimate:
\begin{align}\label{6.1}
\left|\|\nabla u\|_{L^2(D_\delta)}^2-\|\phi_c\|^2\right|\lesssim\omega\left|\langle\phi,\varphi_{0,\delta}\rangle\right|^2.
\end{align}
\end{lemma}
\begin{proof}
Using the divergence theorem in $D_\delta$, we have that
\begin{align*}
\int_{D_\delta}|\nabla u|^2dx&=\overline{k_c}^2\int_{D_\delta}|u|^2dx+\left|\int_{\partial {D_\delta}}u\overline{\frac{\partial
u}{\partial\nu}\Big|_-}d\sigma\right|\\
&=\overline{k_c}^2\int_{D_\delta}|u|^2dx+\int_{\partial D_\delta}\mathcal {S}_{D_\delta}^{k_c}[\phi]\overline{\left(\frac{1}{2}\mathcal{I}+(\mathcal{K}_{D_\delta}^{k_c})^*\right)[\phi]}\, d\sigma.
\end{align*}
Next, by using the same argument as that in the proof of Lemma \ref{lem4.1}, one can derive the estimate (\ref{6.1}).

The proof is complete.
\end{proof}


\begin{theorem}\label{thminter}
Let $u$ be the solution of (\ref{2.8}) inside $D_\delta$. Assume that $|\rho|^{-1}\omega^2\delta\leq c_1$ for a sufficiently small $c_1$, under conditions (C1) and (C2), we have the following results.

\begin{enumerate}
\item[(1)] If $\frac{1}{2}\frac{\theta+\varepsilon_m^{-1}}{\theta-\varepsilon_m^{-1}}\neq\lambda_{j,\delta}$ for any $j\geq0$, then there exists a constant $C$ independent of $\delta$ such that
\begin{equation}
\left\|\nabla u\right\|_{L^2(D_\delta)}\leq C.
\end{equation}

\item[(2)] If $\frac{1}{2}\frac{\theta+\varepsilon_m^{-1}}{\theta-\varepsilon_m^{-1}}=\lambda_{j,\delta}$ and $\langle d\cdot\nu,\varphi_{j,\delta}\rangle\neq0$ for some $j\geq1$,
it follows that
\begin{equation}\label{}
\left\|\nabla u\right\|_{L^2(D_\delta)}\gtrsim\mathcal{O}(|\rho|^{-1}\omega)+o(|\rho|^{-1}\omega\delta)
+\mathcal{O}(|\rho|^{-1}\omega^2)+\mathcal{O}(\omega^{\frac{1}{2}}).
\end{equation}
Furthermore, assuming that $|\rho|=o(\omega)$ and $|\rho|^{-1}\omega^2\delta\rightarrow 0$ as $\omega\rightarrow 0$, then it holds that
\begin{equation}\label{}
\left\|\nabla u\right\|_{L^2(D_\delta)}\rightarrow\infty\quad \mbox{as}\ \ \omega\rightarrow 0.
\end{equation}
\end{enumerate}
\end{theorem}

\begin{proof}
(1) From the relation (\ref{relation}), one sees that $\|\phi\|\leqslant\|\psi\|+\|\varphi_{0,\delta}\|+\mathcal{O}(\omega)$. Noticing the estimate of $\psi$, we have
\[
\left\|\nabla u\right\|_{L^2(D_\delta)}^2=\left\|\nabla \mathcal {S}_{D_\delta}^{k_c}[\phi]\right\|_{L^2(D_\delta)}^2\lesssim\left\|\phi\right\|^2\lesssim\|\psi\|^2+\|\varphi_{0,\delta}\|^2+\mathcal{O}(\omega^2)\lesssim C.
\]

(2) Since
\[
\phi_c+\langle\phi,\varphi_{0,\delta}\rangle\varphi_{0,\delta}=\phi=\psi+c\varphi_{0,\delta}+\mathcal{O}(\omega)
=\psi_c+(c+\langle\psi,\varphi_{0,\delta}\rangle)\varphi_{0,\delta}+\mathcal{O}(\omega),
\]
one can deduce that
\begin{align*}
\langle\phi,\varphi_{0,\delta}\rangle&=(c+\langle\psi,\varphi_{0,\delta}\rangle)\varphi_{0,\delta}+\mathcal{O}(\omega),\\
\langle\phi_c,\varphi_{j,\delta}\rangle&=\langle\psi_c,\varphi_{j,\delta}\rangle+\mathcal{O}(\omega).
\end{align*}
By Lemma \ref{interle01}, we have that
\begin{align*}
\left\|\nabla u\right\|_{L^2(D_\delta)}^2&\gtrsim
\left\|\phi_c\right\|^2-\omega\left|\langle\phi,\varphi_{0,\delta}\rangle\right|^2\\
&\gtrsim\left\|\psi_c+(c+\langle\psi-\phi,\varphi_{0,\delta}\rangle)\varphi_{0,\delta}+\mathcal{O}(\omega)\right\|^2-\omega\left|\langle\phi,\varphi_{0,\delta}\rangle\right|^2\\
&\gtrsim\left(\left\|\psi_c\right\|-(c+\langle\psi-\phi,\varphi_{0,\delta}\rangle)\left\|\varphi_{0,\delta}\right\|-\mathcal{O}(\omega)\right)^2
-\omega\left(\left|c+\langle\psi,\varphi_{0,\delta}\rangle\right|+\mathcal{O}(\omega)\right)^2\\
&\gtrsim\left\|\psi_c\right\|^2-\omega\left|\langle\psi,\varphi_{0,\delta}\rangle\right|^2+\mathcal{O}(\omega)\\
&\gtrsim|\rho|^{-2}\left|\langle f,\varphi_{j,\delta}\rangle\right|^2-\omega\left|\langle \psi,\varphi_{0,\delta}\rangle\right|^2+\mathcal{O}(\omega).
\end{align*}
Hence, following the proof of (2) of Theorem \ref{thm4.3}, one can prove the assertion.

\end{proof}

\begin{remark}\label{remee}
Note that $|\rho|=\Im(\varepsilon_c)/|\varepsilon_c|^2$ and recall the definition of the electric energy $\mathcal{E}$ in (\ref{ee}).
{From (2) of Theorem \ref{thminter}, we have}
\begin{equation}
\mathcal{E}\gtrsim\mathcal{O}(|\rho|^{-1}\omega^2)+o(|\rho|^{-1}\omega^2\delta^2)+\mathcal{O}(|\rho|^{-1}\omega^4)
+\mathcal{O}(|\rho|\omega^2).
\end{equation}
Hence, if we choose that {$|\rho|^{-1}\omega^2\delta\rightarrow0$}, when $\frac{1}{2}\frac{\theta+\varepsilon_m^{-1}}{\theta-\varepsilon_m^{-1}}=\lambda_{j,\delta}$ for some $j\geq1$, and furthermore $\langle d\cdot\nu,\varphi_{j,\delta}\rangle\neq0$ and $|\rho|=o(\omega^2)$ (e.g. one can select $|\rho|=\omega^3$, $\delta=\omega^\varrho$, ($\varrho>1$)), then the electric energy $\mathcal{E}$ blows up (which is dissipated into heat). Moreover, at the same time, by Theorems \ref{thm4.3} and \ref{thminter}, we see that the scattering energy (namely $\|\nabla u^s\|_{L^2(\mathbb{R}^3\backslash\overline{D_\delta)}}$) and the internal energy (namely $\|\nabla u\|_{L^2(D_\delta)}$) also blow up. Hence, we actually have that when resonance occurs, all of the three resonance criteria in Definition~\ref{depr} can be fulfilled.
\end{remark}

\subsection{Simple numerical illustrations}\label{sect:4.3}

\begin{figure}[htp]
\begin{center}
{\includegraphics[width=0.3\textwidth]{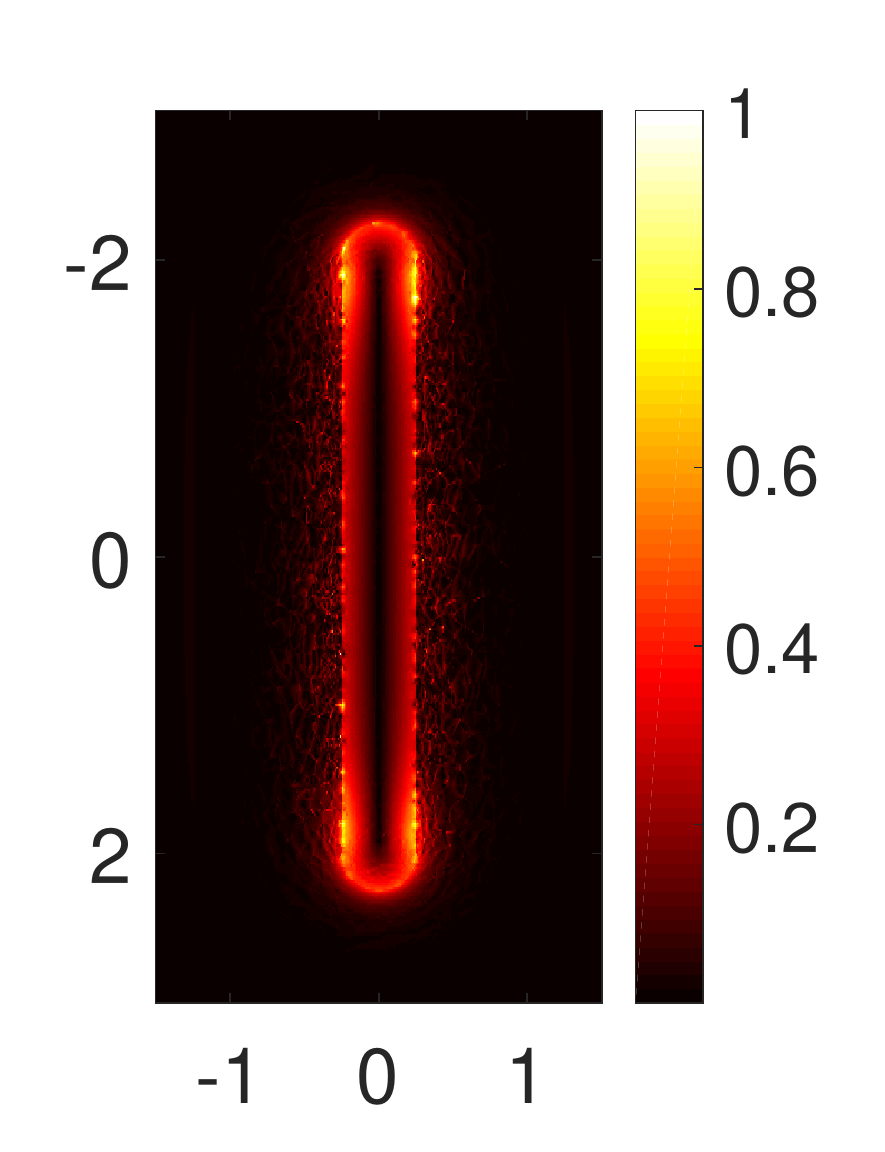}}
{\includegraphics[width=0.3\textwidth]{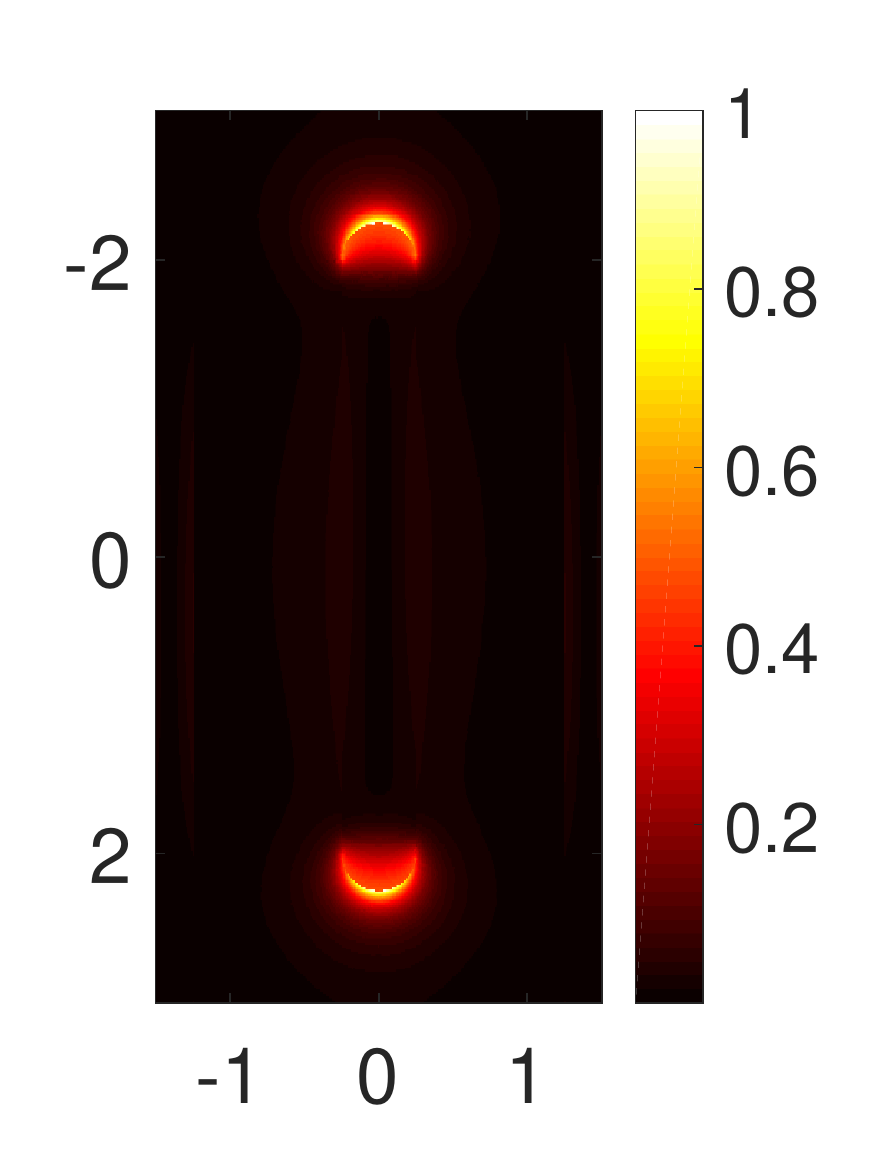}}
\end{center}
\caption{\label{fig:Plasmon03} Left: Slice plot of the normalized resonance strength for the scattering field in the middle figure of Fig.~\ref{fig:Plasmon01} on the $(x_2,x_3)$-plane; $\mathbb{E}_o\approx 1.19\times 10^2$ and  $\mathbb{E}_i\approx 1.88\times 10^2$; Right: Slice plot of the normalized resonance strength for the scattering field in the right figure of Fig.~\ref{fig:Plasmon01} on the $(x_2,x_3)$-plane; $\mathbb{E}_o\approx 1.69\times 10^3$ and  $\mathbb{E}_i\approx 8.44\times 10^2$.}
\end{figure}

\begin{figure}[htp]
\begin{center}
{\includegraphics[width=0.3\textwidth]{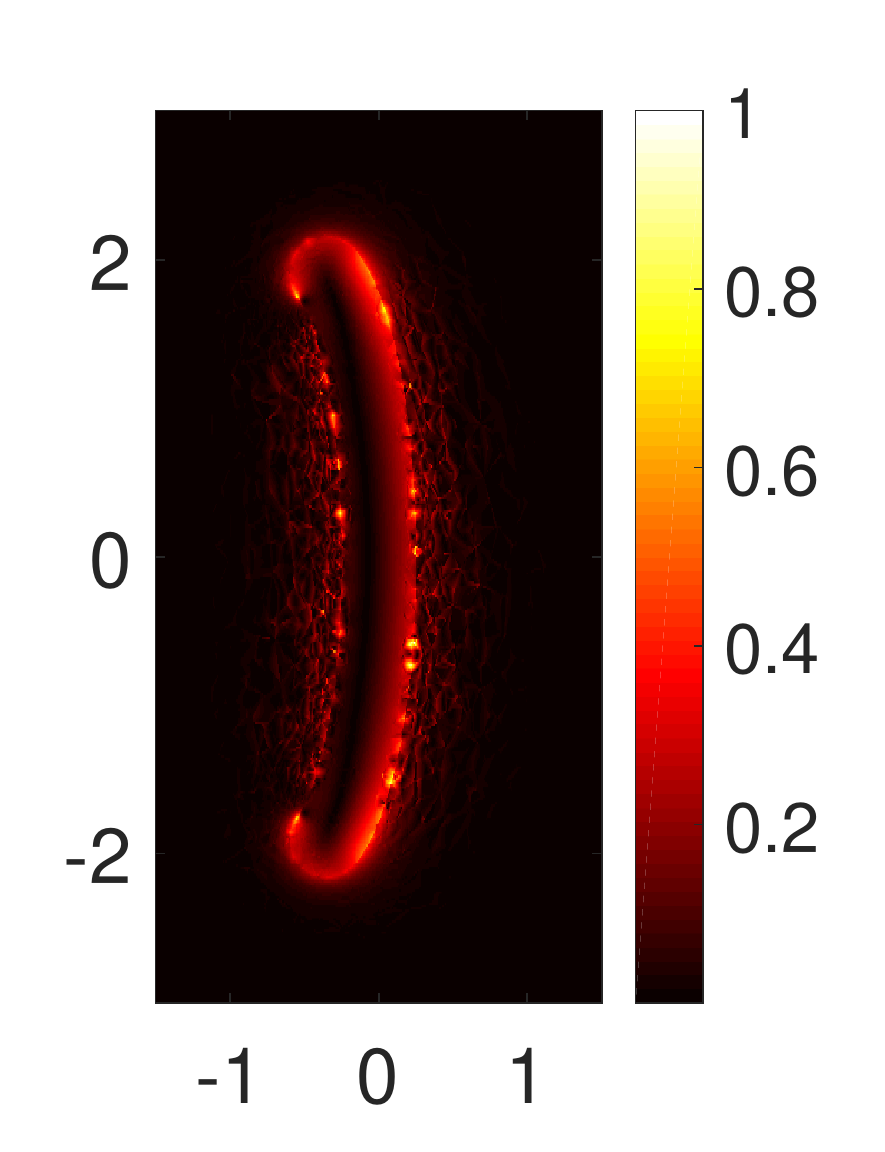}}
{\includegraphics[width=0.3\textwidth]{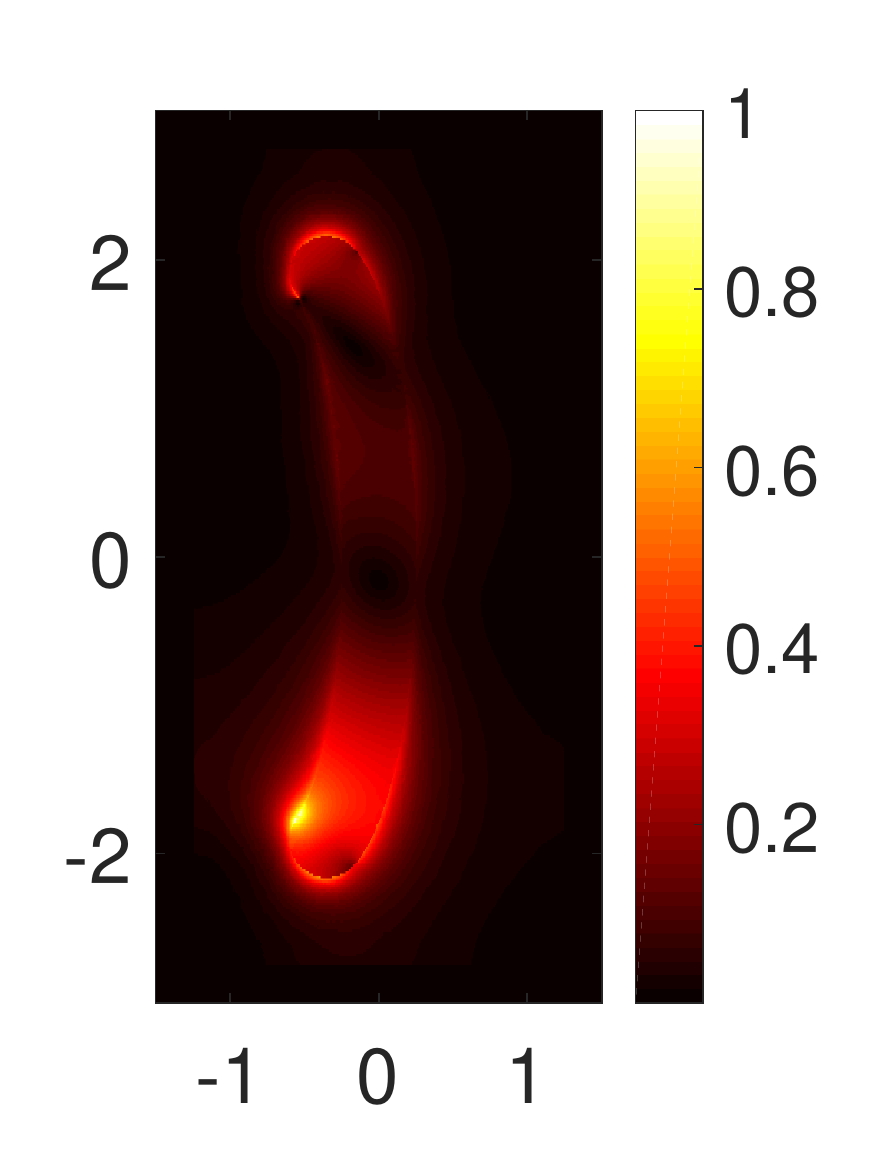}}
\end{center}
\caption{\label{fig:Plasmon04} Left: Slice plot of the normalized resonance strength for the scattering field in middle figure of Fig.~\ref{fig:Plasmon02} on the $(x_2,x_3)$-plane; $\mathbb{E}_o\approx 1.97\times 10^2$ and  $\mathbb{E}_i\approx 1.80\times 10^2$;  Right: Slice plot of the normalized resonance strength for the scattering field in middle figure of Fig.~\ref{fig:Plasmon02} on the $(x_2,x_3)$-plane; $\mathbb{E}_o\approx 1.87\times 10^3$ and  $\mathbb{E}_i\approx 1.91\times 10^3$.
}
\end{figure}

Up till now, we have analyzed the plasmon resonance of the Helmholtz system \eqref{2.8} associate with the nanorod $D_\delta$ in terms of the blowup of the external energy (Theorem~\ref{thm4.3}), the internal energy (Theorem~\ref{thminter}) and the electric energy (Remark~\ref{remee}). Following similar analysis to that in Section~\ref{sect:2.4}, one can show more quantitative properties of the resonant fields, $\nabla u^s(x)$ for $x\in\mathbb{R}^3\backslash\overline{D_{\delta}}$ and $\nabla u(x)$ for $x\in D_\delta$. In particular, the anisotropy of the resonant behaviours due to the anisotropic geometry of $D_\delta$. In principle, one can expect that the resonance strength, namely $|\nabla u^s(x)|$ should be stronger in the vicinity of the two end-parts than that in the vicinity of the facade-part of the nanorod $D_\delta$. To achieve a concrete and vivid idea of such anisotropic resonant behaviours, under the same setup of the numerical examples in Figures \ref{fig:Plasmon01} and \ref{fig:Plasmon02}, we numerically compute the corresponding gradient fields and present them in Fig.~\ref{fig:Plasmon03} and Fig.~\ref{fig:Plasmon04}. Before that, we notice from the material configuration in \eqref{eq:para2} one has
\[
\frac{1}{2}\frac{\theta+\varepsilon_m^{-1}}{\theta-\varepsilon_m^{-1}}=\frac{1}{2}\cdot\frac{-1+(1+\omega^8)}{-1-(1+\omega^8)}\rightarrow0,\ \ \text{as}\ \omega\rightarrow0.
\]
Owing to $0$ is the unique accumulation point of the eigenvalues of $\mathcal{K}_{D_\delta}^*$, the resonance condition $\frac{1}{2}\frac{\theta+\varepsilon_m^{-1}}{\theta-\varepsilon_m^{-1}}=\lambda_{j,\delta}$ in general does not holds exactly, but holds approximately for $\omega\ll1$. In order to fulfil the resonance condition exactly, one needs to compute the eigenvalues of $\mathcal{K}_{D_\delta}^*$, from which to determine an appropriate $\Re\varepsilon_c$. This is not an easy task and is out of the scope of this paper as we mentioned before that numerical study is not our focus. Nevertheless, according to our analysis made above, we know that there are many values of $\Re\varepsilon_c$ centering around $-1$ such that the resonance condition is fulfilled exactly. We vary $\Re\varepsilon_c$ around $-1$ and the obtained numerical behaviours are similar to the case with $\Re\varepsilon_c=-1$. Hence it is unobjectionable that we stick to the material configuration in \eqref{eq:para2} for our numerical simulations of the resonant fields. With such a clarification, we can proceed with the numerical experiments. Set
\[
\Theta(x)=|\nabla \Re u^s(x)|\chi_{\mathbb{R}^3\backslash\overline{D_\delta}}(x)+|\nabla \Re(u(x)-u^i(x))|\chi_{D_\delta}(x),\ \ x\in\mathbb{R}^3,
\]
which is referred to as the resonance strength of the scattering field. In order to have a better display, similar to the plotting in Fig.~\ref{fig:Plasmon01} and Fig.~\ref{fig:Plasmon02}, we shall normalize $\Theta(x)$ such that its maximum value is 1. Let $\Theta_{\mathrm{normal}}$ be the normalized scattering strength, which is a scalar function and can be used to depict the quantitative resonant behaviours for the scattering wave field.
In addition, we define
\[
\mathbb{E}_o=\|\nabla \Re u^s(x)\|_{L^2(B\setminus{D_\delta})}, \quad \mathbb{E}_i=\|\nabla \Re u(x)\|_{L^2(D_\delta)},
\]
where $B$ is the computing region in our numerical experiment. $\mathbb{E}_o$ and $\mathbb{E}_i$ are respectively the exterior and interior scattering energies of the resonant wave field.
In Fig.~\ref{fig:Plasmon03} and Fig.~\ref{fig:Plasmon04}, we respectively plot the normalized resonant strengths of the wave fields in Fig.~\ref{fig:Plasmon01} and Fig.~\ref{fig:Plasmon02}. It can be easily seen that the resonant strength is generically stronger in the vicinity of the two end-parts than that in the vicinity of the facade-part of the nanorod.



\section{Conclusion}

In this paper, we considered the plasmon resonance associated with the Hemholtz system in the quasi-static regime. The plasmon resonance is a delicate and subtle resonant phenomenon and its occurrence and quantitative behaviours critically depend on the metamaterial parameters (namely, the negative material parameters), the frequency of the impinging wave field and the geometry of the nanostructures that are coupled together in a highly intricate manner. Most of the existing theoretical studies are concerned with nanostructures of isotropic geometries that are uniformly small in all dimensions, which help to diminish the geometric influence in the resonance analysis. We investigated the plasmon resonance for curved nanorods, which present anisotropic geometries. The anisotropic geometries present significant challenges for the resonance analysis. By developing novel asymptotic analysis and spectral analysis techniques, we established comprehensive resonance and non-resonance results for this type of nanostructures through carefully analyzing the wave fields and their quantitative behaviours inside and outside the nanostructures. Moreover, through our delicate and subtle analysis, we could quantitatively characterize the anisotropy of the resonant behaviours that is caused by the anisotropic geometries of the curved nanorods. It is mentioned that the geometric properties of the plasmon resonance were also studied in two recent articles \cite{ACL20,BLLW}. In fact, it is shown that the Neumann-Poincar\'e (NP) eigenfunctions, namely the eigenfunctions for the Neumann-Poincar\'e operator $K_{\Sigma}^*$, localize at high-curvature places of $\partial\Sigma$ in the high-mode-number limit. The localizing properties of the NP eigenfunctions indicate that the plasmon resonance is stronger at places on $\partial\Sigma$ where the curvature is sufficiently high, provided the resonance field is caused by the high-mode-number NP eigenfunctions. It is emphasized that the geometric setup considered in \cite{ACL20,BLLW} cannot include the one considered in the current article. Indeed, it is required in \cite{ACL20,BLLW} that the high-curvature place on $\partial\Sigma$ is again of a locally isotropic nature where the curvature is uniformly high in all directions. Finally, we would like to point out that the mathematical methods developed in this work pave the way for many subsequent developments in studying the quantitative geometric properties of the plasmon resonances in different setups as well as considering the corresponding applications, in particular the one for the full Maxwell system that governs the general electromagnetic wave scattering.

\section*{Acknowledgments}

The work of Y. Deng was supported by NSF grant of China No. 11971487 and NSF grant of Hunan No. 2020JJ2038.
The work of H. Liu was supported by a startup fund from City University of Hong Kong and the Hong Kong RGC General Research Funds, 12301218, 12302919 and 12301420.
The work of G. Zheng was supported by NSF grant of China No. 11301168. The authors would like to acknowledge the great helps from Dr. Xianchao Wang on the numerical simulations in this work.


\end{document}